\documentclass{amsart}
\usepackage{graphicx}
\usepackage[]{epsfig}
\usepackage{amssymb}
\usepackage[]{pstricks}
\newpsobject{malla}{psgrid}{subgriddiv=1,griddots=10,gridlabels=6pt}
\newcommand{\lan}{\langle}
\newcommand{\ran}{\rangle}
\newtheorem{theorem}{Theorem}[section]
\newtheorem{proposition}[theorem]{Proposition}
\newtheorem{corollary}[theorem]{Corollary}
\newtheorem{lemma}[theorem]{Lemma}
\newtheorem{remark}[theorem]{Remark}

\newcommand{\p}{\partial}

\newcommand{\R}{{\mathbb R}}
\newcommand{\Z}{{\mathbb Z}}

\newcommand{\T}{{\mathbb T}}
\title[Well-Posedness for the 1-D Schr\"odinger-Debye System ]
{Sharp bilinear estimates and well-posedness for the 1-D
Schr\"odinger-Debye system}

\author{Ad\'an J. Corcho}
\author{Carlos Matheus}
\address{\textbf{Ad\'an J. Corcho}
\newline
Universidade Federal de Alagoas.
\newline
Instituto de Matem\'atica.
\newline
Campus A. C. Sim\~oes, Tabuleiro dos Martins, 57072-900.
\newline
Macei\'o-AL-Brazil.} \email{adan@mat.ufal.br}

\address{\textbf{Carlos Matheus}
\newline
Coll\`ege de France, 3, Rue d'Ulm, CEDEX 05.
\newline
Paris, France.} \email{matheus@impa.br}

\thanks{This research was supported by CAPES, Brazil and ANR, France.}


\date{November 8, 2008.}

\begin{document}
\maketitle

\setcounter{page}{1}

\begin{quote}
{\normalfont\fontsize{8}{10}\selectfont {\bfseries Abstract.} We
establish local and global well-posedness for the initial value
problem associated to the one-dimensional Schr\"odinger-Debye (SD)
system for data in the Sobolev spaces with low regularity. To obtain local
results we prove two new sharp bilinear estimates for the coupling terms of
this system in the continuous and periodic cases. Concerning global results,
in the continuous case, the system is shown to be globally well-posed in
$H^s\times H^s, -3/14< s< 0$. For initial data in Sobolev spaces
with high regularity ($H^s\times H^s,\; s>5/2$),
Bid\'egaray \cite{Bidegaray} proved that there are one-parameter families of solutions of the
SD system converging to certain solutions of the cubic \emph{nonlinear  Schr\"odinger
equation} (NLS). Our results bellow $L^2\times L^2$ say that
the SD system is not a good approach of  the cubic NLS in Sobolev  spaces with low
regularity, since the cubic NLS is known to be ill-posed below $L^2$.
The proof of our global result uses the \textbf{I}-method introduced by
Colliander, Keel, Staffilani, Takaoka and Tao. \par}
\end{quote}

\section{\textbf{Introduction}}

This paper is devoted to the Initial Value Problem(IVP) for the
Schr\"odinger-Debye system, that is,
\begin{equation}
\label{S-Debye}
\begin{cases}
i\p_tu+ \tfrac{1}{2}\partial_x^2u=uv, \quad t  \in {\mathbb R},\quad x\in M\\
\sigma \p_tv + v = \epsilon |u|^2,\\
u(x,0)=u_0(x) , \quad v(x,0)=v_0(x),
\end{cases}
\end{equation}
where $u=u(x,t)$ is a complex valued function, $v=v(x,t)$ is a
real valued function,  $\sigma >0$, $ \epsilon =\pm 1$ and $M$ is
the real line $\R$ (continuous case) or the torus $\T$ (periodic
case).

The  well-posedness for the IVP (\ref{S-Debye}) with initial data in the classical
Sobolev spaces $H^{k}(M) \times H^{s}(M)$  was studied recently by Corcho and Linares \cite{Corcho} when $M=\R^n (n=1,2,3)$ and by  Arbieto and Matheus \cite{AM} when $M=\T^n$. Specifically, in the one-dimensional case they obtained the following results:
\begin{itemize}
    \item \small{local well-posedness in  $H^{s}(\R) \times H^{s}(\R)$ for $0< s <
    1$;}

    \item  \small{global well-posedness in $L^2(\R) \times
    L^2(\R)$ and  $H^{\frac{1}{2}}(\R) \times
    L^2(\R)$;}

    \item  \small{global well-posedness in $H^{k}(\R) \times H^{s}(\R)$ for $k-1/2< s
    \leq
    k$\;\;and\;\;$1/2< k \le 1$;}

    \item \small{local and global well-posedness in $H^{s}(\T) \times H^{s}(\T)$\; for $s\ge
    0$}.
\end{itemize}

The proof of the theorems in the works \cite {Corcho} and \cite{AM} uses Picard fixed-point method in certain spaces. To do so, the authors start by \emph{decoupling} the SD system (\ref{S-Debye}), i.e., they write:
\begin{equation}
u(t) = U(t)u_0 -i\int_0^t U(t-t')\left(e^{-\frac{t'}{\sigma}} v_0 u(t') + \frac{\varepsilon}{\sigma} u(t') \int_0^{t'} e^{-\frac{(t'-t'')}{\sigma}} |u(t'')|^2 dt''\right) dt',
\end{equation}
where $U(t)= e^{it\Delta/2}$ is the Schr\"odinger linear unitary
group. In the sequel, they prove some multilinear estimates for the
nonlinearities in order to Picard's argument run correctly, i.e.,
they show a \emph{bilinear} estimate for the term
$$ \int_0^t U(t-t')\cdot e^{-\frac{t'}{\sigma}} v_0 u(t') dt'$$
and a \emph{trilinear} estimate for the term
$$\int_0^t U(t-t')\frac{\varepsilon}{\sigma} u(t') \left(\int_0^{t'} e^{-\frac{(t'-t'')}{\sigma}} |u(t'')|^2 dt''\right) dt'.$$

Analogously to~\cite{Corcho} and ~\cite{AM}, we are interested in the
local well-posedness of IVP (\ref{S-Debye}) for initial data with
low regularity for $M=\T$ and $M=\R$, specially local and global
well-posedness in the continuous case with initial data in $H^k\times
H^s$ for negative Sobolev indices $(k,s)$. Unfortunately, it is not
reasonable to expect that the approach discussed above can be pushed
to work with negative Sobolev indices. Indeed, similarly to the
situation of Schr\"odinger (NLS) equation, we know that such
trilinear estimates holds only for non-negative indices.

Bearing the difficulty in mind, we propose in this paper a slightly different approach: instead of decoupling the SD system before studying its integral formulation (which leads to trilinear estimates), we keep the SD system coupled so that we have only to deal with bilinear estimates (for the coupling terms $uv$ and $|u|^2$). To understand what is the advantage of our new proposal, we review the bilinear estimates for the quadratic NLS obtained by Kenig, Ponce and Vega.

In \cite{KPV} Kenig, Ponce and Vega considered the initial value
problem
\begin{equation}
\label{NLS}
\begin{cases}
i\p_tu+ \p_x^2u=\alpha N_j(u,u), \quad x,t\in \R, \quad j=1,2,3\\
u(x,0)=u_0(x),
\end{cases}
\end{equation}
where $N_1(u,u)=u\bar u$, $N_2(u,u)=\bar {u}^2$ and
$N_3(u,u)=u^2$. They established the following sharp bilinear
estimates:

\begin{enumerate}
\item[($\textbf{B}_1$)] $\|N_1(u,u)\|_{X^{s,b-1}}\lesssim
\|u\|^2_{X^{s,b}},\;\;\text{for}\;\; s > -1/4 \;\;\text{and}\;\;
b>1/2$;

\vspace{0.2cm}
\item [($\textbf{B}_2$)] $\|N_j(u,u)\|_{X^{s,b-1}}\lesssim \|u\|^2_{X^{s,b}},\;\;
\text{for}\;\; s > -3/4\;\;\text{and}\;\; b>1/2 ,\;\text{with}\;\;
j=2,3,$
\end{enumerate}
where
\begin{equation}\begin{split}\label{Bourgain-Space}
\|f\|_{X^{s,b}}&=\|U(-t)f\|_{H^b_t(\R,H^{s}_x)}\\
&=\left(\int_{\R^2} (1+ |\xi|)^{2s}(1+ |\tau+
\xi^2|)^{2b}|\widehat{f}(\xi,\tau)|^2d\xi d\tau \right)^{1/2}
\end{split}\end{equation} and $U(t):=e^{it\partial_x^2}$ is the
corresponding Schr\"odinger generator (unitary group) associated
to the linear problem. Using the estimates ($\textbf{B}_1$) and
($\textbf{B}_2$) and properties of the $X^{s,b}$ spaces together
with the contraction mapping principle they proved local
well-posedness for (\ref{NLS}) in $H^s(\R)$ for $s> -1/4$ ($j=1$)
and for $s> -3/4$ ($j=2,3$).

Similar results were given  in the periodic case, where $\|\cdot
\|_{X^{s,b}_{per}}$ is defined by
\begin{equation}\begin{split}\label{Bourgain-Space-Periodic}
\|f\|_{X^{s,b}_{per}}&=\left(\sum_{n\in
\Z}\int_{-\infty}^{+\infty} (1+ |n|)^{2s}(1+ |\tau+
n^2|)^{2b}|\widehat{f}(n,\tau)|^2d\tau \right)^{1/2}
\end{split}\end{equation}
and the corresponding bilinear estimates obtained are the
followings:
\begin{enumerate}
\item[($\textbf{B}_3$)] $\|N_1(u,u)\|_{X^{s,b-1}_{per}}\lesssim
\|u\|^2_{X^{s,b}_{per}},\;\;\text{for}\;\; s \ge 0
\;\;\text{and}\;\; b\in(1/2,1)$;

\vspace{0.2cm}
\item [($\textbf{B}_4$)] $\|N_j(u,u)\|_{X^{s,b-1}_{per}}\lesssim \|u\|^2_{X^{s,b}_{per}},\;\;
\text{for}\;s > -1/2\;\; \text{and}\; b\in (1/2,1)
,\;\text{with}\; j=2,3.$
\end{enumerate}

As explained above, in our case the nonlinear interactions are $uv$ and $|u|^2$. These terms are similar to $N_3$ and $N_1$, respectively,
but the characteristics of linear part of each equation involved
in the system (\ref{S-Debye}) are antisymmetric. Therefore, our task is to find new mixed bilinear
estimates for the coupling terms $uv$ and $|u|^2$.

Before stating the results we will give some useful notations. Let
$\psi$ be a function in $C_0^{\infty}$ such that $0\le \psi(t) \le
1$,
$$\psi(t)=
\begin{cases}
1&\;\; \text{if}\;\;|t|\le 1,\\
0&\;\; \text{if}\;\;|t|\ge 2,
\end{cases}$$
and $\psi_{T}(t)=\psi (\tfrac{t}{T})$. We denote by $\lambda \pm$
a number slightly larger, respectively smaller, than $\lambda$ and
by $\lan \cdot \ran$ the number $\lan \cdot \ran= 1+|\cdot|$. The
characteristic function on the set $A$ is denoted by $\chi_{A}$.

The next statements show the main local-in-time results achieved in this work.

\begin{theorem}\label{local-theorem-continuous}
For any $(u_0,v_0)\in H^k(\R) \times H^s(\R)$ provided the
conditions:
\begin{equation}\label{local-theorem-continuous-a}
|k|-1/2\le s < \min\{k+1/2, 2k+1/2\}\;\;\; and \;\;\; k> -1/4.
\end{equation}
there exist a positive time $T=T(\|u_0\|_{H^k}, \|v_0\|_{H^s})$ and
a unique solution $(u(t),v(t))$ of the initial value problem
(\ref{S-Debye}) on the time interval $[0,T]$, satisfying
\begin{enumerate}
\item [(i)\;] $\left(\psi_T(t)u, \psi_T(t)v\right)\in X^{k,\tfrac{1}{2}+}\times H^{\tfrac{1}{2}+}(\R, H^s_x)$;
\vspace{0.2cm}
\item [(ii)] $(u, v)\in C\left([0,T]; H^k(\R) \times H^s(\R)\right)$.
\end{enumerate}
Moreover, the map $(u_0,v_0) \longmapsto (u(t),v(t))$ is locally
Lipschitz from $H^k(\R) \times H^s(\R)$ into $C([0,T]; H^k(\R)
\times H^s(\R))$.
\end{theorem}

\begin{theorem}\label{local-theorem-periodic}
For any $(u_0,v_0)\in H^k(\T) \times H^s(\T)$ provided the
conditions:
\begin{equation}\label{local-theorem-periodic-a}
0\le s\le 2k\;\;\; and \;\;\; |s-k|<1.
\end{equation}
there exist a positive time $T=T(\|u_0\|_{H^k}, \|v_0\|_{H^s})$
and a unique solution $(u(t),v(t))$ of the initial value problem
(\ref{S-Debye}), satisfying
\begin{enumerate}
\item [(i)\;] $\left(\psi_T(t)u, \psi_T(t)v\right)\in X_{per}^{k,\tfrac{1}{2}+}
\times H^{\tfrac{1}{2}+}(\R, H^s_{per})$;
\vspace{0.2cm}
\item [(ii)] $(u, v)\in C\left([0,T]; H^k(\T) \times H^s(\T)\right)$.
\end{enumerate}
Moreover, the map $(u_0,v_0) \longmapsto (u(t),v(t))$ is locally
Lipschitz from $H^k(\T) \times H^s(\T)$ into $C([0,T]; H^k(\T)
\times H^s(\T))$.
\end{theorem}

In figures 1 and 2 below, respectively, we design the regions on the $(k,s)$-plane where our local well-posedness theorems in the continuous and periodic settings, respectively, are valid.

Finally, we show that the system (\ref{S-Debye}) is globally
well-posed for a class of data without finite mass, more
precisely:
\begin{theorem}\label{global-theorem}
For any $(u_0,v_0)\in H^s(\R)\times H^s(\R),\; -3/14 < s <0$, the
local solution given in Theorem \ref{local-theorem-continuous} can
be extended to any time interval $[0,T]$ (preserving the
properties (i) and (ii).)
\end{theorem}

The plan of this paper is as follows. In Section
\ref{section-preliminary} are given preliminary estimates needed
to establish the new mixed bilinear estimates for coupling terms
of system (\ref{S-Debye}) and the proof of these estimates will be
given in Sections \ref{section-bilinear-continuous} and
\ref{section-bilinear-periodic}. Moreover, we observe that our
local results, given in theorems \ref{local-theorem-continuous}
and \ref{local-theorem-periodic}, are consequences of these
bilinear estimates by using the standard contraction mapping
principle and the properties of $X^{s,b}$ spaces. For instance, see the works \cite{KPV}, \cite{Bekiranov1} and \cite{Ginibre}. Finally, in Section \ref{section-global-results} we proof Theorem
\ref{global-theorem} using the I-method combined with the
following \emph{refined} Strichartz type estimate for the Schr\"odinger equation:
\begin{equation}\label{Strichartz-Estimate}
\|(D_x^{1/2}f)\cdot g\|_{L_{xt}^2}\lesssim \|f\|_{X^{0,1/2+}} \|g\|_{X^{0,1/2+}},
\end{equation}
if $|\xi_1|\gg |\xi_2|$ for any $|\xi_1|\in \textrm{supp}(\widehat{f}), |\xi_2|\in \textrm{supp}(\widehat{g})$. See~\cite{CKSTT3} and~\cite{Grunrock} for more details about refined Strichartz estimates.

We finish with the following interesting remark: in the work
\cite{Bidegaray} it was shown that as the parameter $\sigma$ tends
to zero, solutions the system (\ref{S-Debye}) converge (in
$H^s(\R)$\; for\; $s> 5/2$) to those of the cubic nonlinear
Schr\"odinger equation. Our local results in Theorem
\ref{local-theorem-continuous} show that this fact is not true in
Sobolev spaces with low regularity since the cubic Schr\"odinger
equation is not locally well-posed below $L^2$ in the continuous
case (in the sense that the associated flow is not uniformly continuous).

\vspace{0.5cm}
\begin{figure}[ht]\label{Region-Continuous}
\centering \psset{unit=1.0cm}
\begin{pspicture}(-2,-2)(4,4)%
\malla%
\psline{->}(-2,0)(4,0)%
\psline{->}(0,-2)(0,4)%
\psline[](-0.25,-0.25)(0,-0.5)(4,3.5)%
\pspolygon[fillstyle=solid,fillcolor=yellow,linewidth=0.5pt]
(-0.25,-0.25)(-0.25,0)(0,0.5)(3.5,4)(4,4)(4,3.5)(0,-0.5)(-0.25,-0.25)%
\psline[linecolor=gray](-0.25,0)(0.5,0)%
\psline[linecolor=gray](0,-0.5)(0,0.5)%
\psline[linecolor=red](2.2,2.2)(4,4)%
\psline[linecolor=red](1.7,1.7)(-0.25,-0.25)%
\rput(2,1.85){$\mathcal{W}$}%
\rput(1,0.77){$\ell$}%
\rput(3.7,-0.25){$k$}%
\rput(0.25,3.8){$s$}%
\rput(-0.3,0.25){\small{$r_3$}}%
\rput(1.7,2.4){\small{$r_2$}}%
\rput(0.3,-0.5){\small{$r_1$}}%
\end{pspicture}
\vspace{0.5cm} \caption{Well-posedness
results for Schr\"odinger-Debye  system in the continuous case
($M=\R$). The region $\mathcal{W}$, limited by the lines $r_1:
|k|-s=1/2$\;and\;  $r_2: s-k=1/2$, $r_3: s-2k=1/2$, for $k\geq -1/4$, contain the indices $(k,s)$ where the local well-posedness is achieved in Theorem \ref{local-theorem-continuous}. Global results, given in Theorem \ref{global-theorem}, are obtained on the line $\ell: s=k$\; for
$-3/14< k\le0$.}
\end{figure}
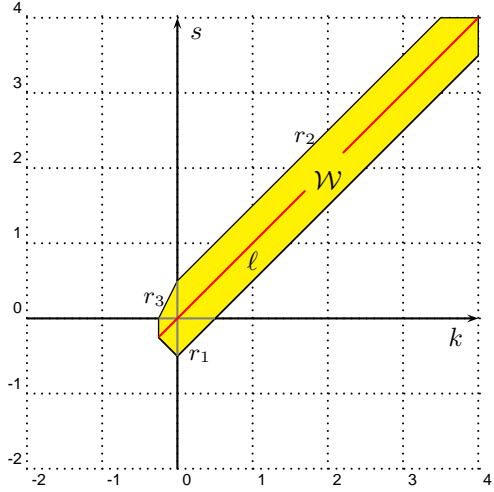

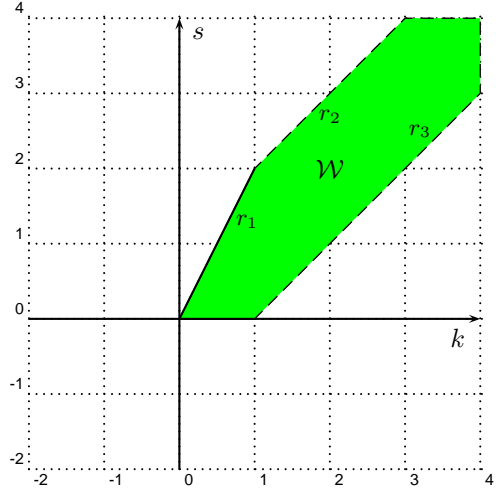
\begin{figure}[h]\label{Region-Periodic}
\vspace{0.7cm} \centering \psset{unit=1.0cm}
\begin{pspicture}(-2,-2)(4,4)%
\malla%
\psline{->}(-2,0)(4,0)%
\psline{->}(0,-2)(0,4)%
\pspolygon[fillstyle=solid,fillcolor=green,linewidth=0.5pt,linestyle=dashed]
(0,0)(1,0)(4,3)(4,4)(3,4)(1,2)(0,0)%
\psline[](0,0)(1,2)%
\psline[](0,0)(1,0)%
\rput(2,2){$\mathcal{W}$}%
\rput(3.7,-0.25){$k$}%
\rput(0.25,3.8){$s$}%
\rput(0.9,1.3){\small{$r_1$}}%
\rput(2,2.7){\small{$r_2$}}%
\rput(3.2,2.5){\small{$r_3$}}%
\end{pspicture}
\vspace{0.5cm} \caption{Well-posedness results for periodic
Schr\"odinger-Debye system ($M=\T$). The region $\mathcal{W}$,
limited for the lines $r_1: s=2k$, $r_2: s=k+1$ and $r_3: s=k-1$,
contain the indices $(k,s)$ where the local well-posedness is
achieved in Theorem \ref{local-theorem-periodic}.}
\end{figure}

\newpage

\section{\textbf{Preliminary Estimates}}\label{section-preliminary}

Firstly, we recall some estimates contained in the
work~\cite{Ginibre} of Ginibre, Tsutsumi and Velo concerning the
Zakharov system:

\begin{lemma}\label{pe-sd-lemma1}
Let $-1/2< b'\leq 0 \leq b\leq b'+1$ and $T\in[0,1]$. Then, for
$F\in H^{b'}_t(\R, H^s_x)$ we have
\begin{equation}\label{pe-sd-lemma1-a}
\|\psi _1(t)\omega_0 \|_{H^b_t(\R, H^s_x)}\leq C\|\omega
_0\|_{H^s},
\end{equation}
\begin{equation}\label{pe-sd-lemma1-b}
\left \|\psi _T(t)\int_0^tF(t',\cdot)dt' \right \|_{H^b_t(\R,
H^s_x)}\leq CT^{1-b+b'}\|F\|_{H^{b'}_t(\R, H^s_x)}.
\end{equation}
\end{lemma}
\begin{proof}
See Lemma 2.1 in \cite{Ginibre}.
\end{proof}

\begin{lemma}\label{l.ginibre}It holds
\begin{equation}
\int_{\mathbb{R}^4}
\frac{|\widehat{f}(\xi,\tau)\widehat{g}(\xi_1,\tau_1)\widehat{h}(\xi_2,\tau_2)|
\langle\xi\rangle^{1/2}}{\langle\sigma\rangle^{d}\langle\sigma_1\rangle^{d_1}\langle\sigma_2\rangle^{d_2}}
d\xi_1 d\tau_1 d\xi d\tau \lesssim
\|f\|_{L^2_{xt}}\|g\|_{L^2_{xt}}\|h\|_{L^2_{xt}},
\end{equation}
where $\xi=\xi_1+\xi_2$, $\tau=\tau_1+\tau_2$, $\sigma:=\tau$,
$\sigma_1:=\tau_1-\frac{1}{2}\xi_1^2$,
$\sigma_2:=\tau_2+\frac{1}{2}\xi_2^2$ and $d,d_1,d_2>1/4$,
$d+d_1>3/4$, $d+d_2>3/4$.
\end{lemma}

\begin{proof}See~\cite[p.422--424]{Ginibre}.
\end{proof}

Next, we recall some elementary calculus inequalities:

\begin{lemma}\label{pe-sd-lemma2} Let $p, q > 0$. Then
for $r=\min\{p,q\}$ with $p+q >1+r$ there exists $C>0$ such that
\begin{equation}\label{pe-sd-lemma2-a}
\int_{-\infty}^{\infty}\frac{dx}{\lan x-\alpha \ran^p \lan x-\beta
\ran^q}\leq \frac{C}{\lan \alpha-\beta \ran^r}.
\end{equation}
Furthermore, for $p>1$ and $q>1/2$ there exists a $C>0$ such that

\begin{equation}\label{pe-sd-lemma2-b}
\int_{-\infty}^{\infty}\frac{dx}{\lan \alpha x - \beta \ran^p}\leq
\frac{C}{|\alpha|},\quad \text{for}\quad \alpha \neq 0,
\end{equation}
\begin{equation}\label{pe-sd-lemma2-c}
\int_{-\infty}^{\infty}\frac{dx}{\lan \alpha_0 + \alpha_1x +
\tfrac{1}{2}x^2\ran^q}\leq C.
\end{equation}
\end{lemma}

\begin{proof}
See the work \cite {Bekiranov1}.
\end{proof}

Finally, we recall some time localization properties of the Bourgain spaces:

\begin{lemma}\label{l.time}Let $-1/2<b'\leq b<1/2$, $s\in\mathbb{R}$ and $0<T<1$. It holds 
$$\|\psi_T(t)u\|_{X^{s,b'}}\lesssim T^{b-b'}\|u\|_{X^{s,b}}$$
and 
$$\|\psi_T(t) v\|_{H_t^{b'}H_x^s}\lesssim T^{b-b'} \|v\|_{H_t^b H_x^s}.$$
\end{lemma}

\begin{proof}
See lemma 2.11 of the book~\cite{T}.
\end{proof}

\section{\textbf{Bilinear Estimates for the Coupling Terms in the Continuous Case}}
\label{section-bilinear-continuous} The aim of this section is the
study of the crucial sharp bilinear estimates for the coupling
terms in the continuous cases. In order to do so, this section is organized as follows: first, we present the proof of the relevant bilinear estimates assuming certain restrictions on the Sobolev indices $s$ and $k$ of the initial data; after this, we show a series of counter-examples showing that our restrictions on $s$ and $k$ are necessary.

\subsection{Proof of the bilinear estimates I: the continuous case}

\begin{proposition}\label{proposition-uv-continuous}
Let $1/4<a<1/2$ and $b>1/2$. The bilinear estimate
\begin{equation}\label{Estimativa-Mixta-1-c}
\|uv\|_{X^{k,-a}}\lesssim \|u\|_{X^{k,b}}\|v\|_{H^{b}_t H^s_x}
\end{equation}
holds if\; \small{$|k|- s\le 1/2$}.
\end{proposition}

\begin{proof} We  define
\begin{equation*}
f(\xi, \tau)= \lan\tau + \tfrac{1}{2}\xi^2\ran^b\lan \xi
\ran^{k}\widehat u(\xi,\tau)\;\;\; \text{and}\;\;\;  g(\xi, \tau)=
\lan\tau \ran^b\lan \xi \ran^{s}\widehat v(\xi,\tau).
\end{equation*}

Then, for $u\in X^{k,b}$\; and\; $v\in H^b_tH^s_x$, the $L^2$ duality and the definition (\ref{Bourgain-Space}) show
that (\ref{Estimativa-Mixta-1-c}) is equivalent to prove

\begin{equation}\label{Desigualdade-Equivalente-1}
 W_{f,g}(\varphi)
\lesssim \|f\|_{L^2}\|g\|_{L^2}\|\varphi\|_{L^2},
\end{equation}
for all $\varphi \in L^2(\R^2)$, where
\begin{equation}\label{Desigualdade-Equivalente-1-a}
\small{W_{f,g}(\varphi)=\int_{\R^4} \frac{\lan \xi
\ran^k\bar{\varphi}(\xi, \tau) f(\xi-\xi_1, \tau -\tau_1)g(\xi_1,
\tau_1)}{\lan\tau + \frac{1}{2}\xi^2\ran^a \lan \xi-\xi_1
\ran^k\lan \tau -\tau_1 + \frac{1}{2}(\xi-\xi_1)^2\ran^b\lan \xi_1
\ran^s\lan \tau_1 \ran^b}d\xi_1 d\tau_1d\xi d\tau}
\end{equation}

To estimate $W_{f,g}$ we split $\R^4$ into three regions
$\mathcal{A}_1$, $\mathcal{A}_2$ and  $\mathcal{A}_3$,
\begin{equation*}\begin{split}
&{\mathcal A}_1=\left \{(\xi,\xi_1, \tau,\tau_1)\in {\mathbb R}^4 ;\; |\xi_1|\leq 1 \right\},\\
&{\mathcal A}_2=\left\{(\xi,\xi_1, \tau,\tau_1)\in {\mathbb R}^4
;\; |\xi_1|> 1\;\text{and}\;
|\xi_1-\xi|\geq \tfrac{1}{8}|\xi_1| \right\},\\
&{\mathcal A}_3=\left\{(\xi,\xi_1, \tau,\tau_1)\in {\mathbb
R}^4;\; |\xi_1|> 1 \;\text{and}\; |\tfrac{1}{2}\xi_1-\xi|\geq
\tfrac{1}{8}|\xi_1| \right\}.
\end{split}\end{equation*}
Since
$$\mathcal{S}=\left \{(\xi,\xi_1, \tau,\tau_1)\in {\mathbb R}^4;\; |\xi_1|> 1,\;
|\xi_1-\xi|< \tfrac{1}{8}|\xi_1|\;\text{and}\;
|\tfrac{1}{2}\xi_1-\xi|< \tfrac{1}{8}|\xi_1|\right \}$$ is empty,
we have that $\R^4={\mathcal A}_1\cup {\mathcal A}_2\cup {\mathcal
A}_3$.\;Indeed if $(\xi,\xi_1, \tau,\tau_1)\in {\mathcal
S}$,\;then
$$
\tfrac{1}{2}|\xi_1|=|\xi_1-\xi-(\tfrac{1}{2}\xi_1 -\xi)|\le
|\xi_1-\xi|+|\tfrac{1}{2}\xi_1-\xi|<\tfrac{1}{4}|\xi_1|,
$$
which is a contradiction.

Note that for any point in $\mathcal{A}_3$ we have the following
algebraic inequality
\begin{equation}\label{Algebraic-Relation}
|\tau+\tfrac{1}{2}\xi^2|+|\tau_1|+|\tau-\tau_1+\tfrac{1}{2}(\xi-\xi_1)^2|\ge|\tfrac{1}{2}\xi_1^2-\xi\xi_1|
=|\xi_1||\tfrac{1}{2}\xi_1-\xi|\ge \tfrac{1}{8}|\xi_1|^2,
\end{equation}
and consequently
\begin{equation}\label{Algebraic-Inequality}
\max \left \{|\tau+\tfrac{1}{2}\xi^2|, |\tau_1|,
|\tau-\tau_1+\tfrac{1}{2}(\xi-\xi_1)^2| \right \}\ge
\tfrac{1}{24}|\xi_1|^2.
\end{equation}

Now we separate $\mathcal{A}_3$ into three parts,
\begin{equation*}\begin{split}
&{\mathcal A}_{3,1} =\left\{(\xi,\xi_1, \tau,\tau_1)\in {\mathcal
A}_3;\; |\tau_1|,\;|\tau-\tau_1 + \tfrac{1}{2}(\xi-\xi_1)^2|
\leq |\tau + \tfrac{1}{2}\xi^2| \right \},\\
&{\mathcal A}_{3,2} =\left \{(\xi,\xi_1, \tau,\tau_1)\in {\mathcal
A}_3;\; |\tau-\tau_1 + \tfrac{1}{2}(\xi-\xi_1)^2| ,\;|\tau +
\tfrac{1}{2}\xi^2|
\leq |\tau_1|\right \},\\
&{\mathcal A}_{3,3} =\left \{(\xi,\xi_1, \tau,\tau_1)\in {\mathcal
A}_3;\; |\tau_1| ,\;|\tau + \tfrac{1}{2}\xi^2| \leq |\tau-\tau_1 +
\tfrac{1}{2}(\xi-\xi_1)^2|\right \},
\end{split}\end{equation*}
so that one of the following $|\tau +
\tfrac{1}{2}\xi^2|$,\;$|\tau_1|$\;or\;$|\tau-\tau_1 +
{\tfrac{1}{2}(\xi-\xi_1)}^2|$\; is larger than
$\tfrac{1}{24}|\xi_1|^2$.

We can now define the sets $\Omega_1={\mathcal A}_1\cup {\mathcal
A}_2\cup {\mathcal A}_{3,1}$,\;\;$\Omega_2={\mathcal
A}_{3,2}$\;\;and\;\;$\Omega_3={\mathcal A}_{3,3}$ and it is clear
that $\R^4=\Omega_1\cup \Omega_2\cup \Omega_3$. Then, we
decompose  the integral in $W$ into the followings
$$W(f,g,\varphi)=W_1+W_2+W_3,$$
where
\begin{equation*}\label{Desigualdade-Equivalente-j}
\begin{split} W_j&= \int_{\Omega_j}\!\!\frac{\lan \xi \ran^k\bar{\varphi}(\xi, \tau)  f(\xi-\xi_1, \tau -\tau_1)g(\xi_1, \tau_1)}
{\lan\tau + \frac{1}{2}\xi^2\ran^a \lan \xi_1 \ran^s\lan \tau_1
\ran^b \lan \xi-\xi_1 \ran^k\lan \tau -\tau_1 +
\frac{1}{2}(\xi-\xi_1)^2\ran^b}d\xi_1 d\tau_1d\xi d\tau ,
\end{split}
\end{equation*}
for $j=1,2,3.$

We begin by estimating $W_1$. For this purpose, we integrate over
$\xi_1$ and $\tau_1$ first and then use the Cauchy-Schwarz and
H\"older inequalities and the Fubini's theorem to obtain
\begin{equation}{\begin{split}\label{Estimative-W1} |W_1|^2
& \le \Biggl\|\frac{\lan \xi \ran^{k}}{\lan \tau +
\frac{1}{2}\xi^2 \ran^{a}}\int_{\R^2}\frac{ f(\xi-\xi_1, \tau
-\tau_1)g(\xi_1, \tau_1)}{\lan \xi_1 \ran^s \lan \xi-\xi_1
\ran^k\lan \tau_1 \ran^b \lan \tau -\tau_1 +
\frac{1}{2}(\xi-\xi_1)^2\ran^b}\chi_{\Omega_1}d\xi_1d\tau_1\Biggl
\|^2_{L^2_{\xi,\tau}}\times \\
&\;\;\;\; \times \|\varphi\|^2_{L^2_{\xi,\tau}}\\
&= \int_{\R^2}\!\!\tfrac{\lan \xi \ran^{2k}d\xi d\tau}{\lan \tau +
\frac{1}{2}\xi^2 \ran^{2a}}\left |\int_{\R^2}\frac{ f(\xi-\xi_1,
\tau -\tau_1)g(\xi_1, \tau_1)}{\lan \xi_1 \ran^{2s} \lan \xi-\xi_1
\ran^{2k}\lan \tau_1 \ran^{2b} \lan \tau -\tau_1 +
\frac{1}{2}(\xi-\xi_1)^2\ran^{2b}}\chi_{\Omega_1}d\xi_1d\tau_1\right
|^2\times \\
&\;\;\;\; \times \|\varphi\|^2_{L^2_{\xi,\tau}}\\
& \le \left \| \frac{\lan \xi \ran^{2k}}{\lan \tau +
\frac{1}{2}\xi^2 \ran^{2a}} \int_{\R^2}\frac{\lan \xi_1 \ran^{-2s}
\lan \xi-\xi_1 \ran^{-2k}} {\lan \tau_1 \ran^{2b} \lan \tau
-\tau_1 + \frac{1}{2}(\xi-\xi_1)^2\ran^{2b}}
\chi_{\Omega_1}d\xi_1d\tau_1 \right \|_{L^{\infty}_{\xi,\tau}}\times \\
& \;\;\;\; \times \|f\|^2_{L^2}\|g\|^2_{L^2}\|\varphi\|^2_{L^2}
\end{split}}
\end{equation}

For $W_2$ we put $\tilde{f}(\xi,\tau):=f(-\xi,-\tau)$, integrate
over $\xi$ and $\tau$ first and follow the same steps as above to
get
\begin{equation}\begin{split}\label{Estimative-W2}|W_2|^2 &
\le \left \| \frac{\lan \xi_1 \ran^{-2s}}{\lan\tau_1\ran^{2b}}
\int_{\R^2}\frac{\lan \xi \ran^{2k} \lan \xi-\xi_1 \ran^{-2k}}
{\lan \tau + \tfrac{1}{2}\xi^2 \ran^{2a} \lan \tau -\tau_1 +
\frac{1}{2}(\xi-\xi_1)^2\ran^{2b}}
\chi_{\Omega_2}d\xi d\tau \right \|_{L^{\infty}_{\xi_1,\tau_1}}\times  \\
& \;\;\;\; \times
\|\tilde{f}\|^2_{L^2}\|g\|^2_{L^2}\|\varphi\|^2_{L^2}.
\end{split}
\end{equation}
Note that $\|\tilde{f}\|^2_{L^2}=\|f\|^2_{L^2}.$

Now we use the change of variables $\tau=\tau_1-\tau_2$ and
$\xi=\xi_1-\xi_2$ to transform the region $\Omega_3$ into the set
$\widetilde {\Omega}_3$, that satisfies

$$\widetilde {\Omega}_3 \subseteq \{(\xi_1, \xi_2, \tau_1, \tau_2)\in
\R^4;\;\; \tfrac{1}{8}|\xi_1|^2\le
|\tfrac{1}{2}\xi_1^2-\xi_1\xi_2|\le
3|\tau_2-\tfrac{1}{2}\xi_2^2|\;\;\text{and}\;\; |\xi_1|>1\}$$

Then $W_3$ can be estimated as follows

\begin{equation}\begin{split}\label{Estimative-W3} |W_3|^2
& \le \left \| \frac{\lan \xi_2
\ran^{-2k}}{\lan\tau_2-\tfrac{1}{2}\xi_2^2\ran^{2b}}
\int_{\R^2}\frac{\lan \xi_1 \ran^{-2s} \lan \xi_1-\xi_2 \ran^{2k}}
{\lan \tau_1 \ran^{2b} \lan \tau_1 -\tau_2 +
\frac{1}{2}(\xi_1-\xi_2)^2\ran^{2a}}
\chi_{\widetilde{\Omega}_3}d\xi_1d\tau_1 \right \|_{L^{\infty}_{\xi_2,\tau_2}}\times \\
& \;\;\;\; \times
\|\tilde{f}\|^2_{L^2}\|g\|^2_{L^2}\|\varphi\|^2_{L^2}.
\end{split}
\end{equation}

From estimates (\ref{Estimative-W1}), (\ref{Estimative-W2}) and
(\ref{Estimative-W3}) it suffices to show that the following
expressions are bounded:

\begin{equation}\label{Estimative-W1-Tilde}
\widetilde{W}_1(\xi,\tau):=\frac{\lan \xi \ran^{2k}}{\lan \tau
+ \frac{1}{2}\xi^2 \ran^{2a}} \int_{\R^2}\frac{\lan \xi_1
\ran^{-2s} \lan \xi-\xi_1 \ran^{-2k}} {\lan \tau_1 \ran^{2b} \lan
\tau -\tau_1 + \frac{1}{2}(\xi-\xi_1)^2\ran^{2b}}
\chi_{\Omega_1}d\xi_1d\tau_1,
\end{equation}

\begin{equation}\label{Estimative-W2-Tilde}
\widetilde{W}_2(\xi_1,\tau_1):=\frac{\lan \xi_1
\ran^{-2s}}{\lan\tau_1\ran^{2b}} \int_{\R^2}\frac{\lan \xi
\ran^{2k} \lan \xi-\xi_1 \ran^{-2k}} {\lan \tau +
\tfrac{1}{2}\xi^2 \ran^{2a} \lan \tau -\tau_1 +
\frac{1}{2}(\xi-\xi_1)^2\ran^{2b}} \chi_{\Omega_2}d\xi d\tau,
\end{equation}
and
\begin{equation}\label{Estimative-W3-Tilde}
\widetilde{W}_3(\xi_2,\tau_2):=\frac{\lan \xi_2
\ran^{-2k}}{\lan\tau_2-\tfrac{1}{2}\xi_2^2\ran^{2b}}
\int_{\R^2}\frac{\lan \xi_1 \ran^{-2s} \lan \xi_1-\xi_2 \ran^{2k}}
{\lan \tau_1 \ran^{2b} \lan \tau_1 -\tau_2 +
\frac{1}{2}(\xi_1-\xi_2)^2\ran^{2a}}
\chi_{\widetilde{\Omega}_3}d\xi_1d\tau_1.
\end{equation}

Now using lemma \ref{pe-sd-lemma2}-(\ref{pe-sd-lemma2-a}) and the
inequalities: $\lan \xi\ran^{2k} \le \lan \xi_1\ran^{2|k|}\lan
\xi-\xi_1\ran^{2k}$ and \newline $\lan \xi_1-\xi_2\ran^{2k} \le
\lan \xi_1\ran^{2|k|}\lan \xi_2\ran^{2k}$, for $k\ge 0$, and $\lan
\xi-\xi_1\ran^{-2k} \le \lan \xi_1\ran^{2|k|}\lan \xi\ran^{-2k}$
and  $\lan \xi_2\ran^{-2k} \le \lan \xi_1\ran^{2|k|}\lan
\xi_1-\xi_2\ran^{2k}$, for $k<0$, we have

\begin{equation}\label{Estimative-W1-Tilde-a}
\widetilde{W}_1(\xi,\tau)\le J_1(\xi,\tau):=\frac{1}{\lan \tau +
\frac{1}{2}\xi^2 \ran^{2a}} \int_{-\infty}^{+\infty}\frac{\lan
\xi_1 \ran^{2|k|-2s}} {\lan \tau
+\frac{1}{2}\xi^2+\frac{1}{2}\xi_1^2-\xi\xi_1\ran^{2b}}
\chi_{\Omega_1}d\xi_1,
\end{equation}
\begin{equation}\label{Estimative-W2-Tilde-a}
\widetilde{W}_2(\xi_1,\tau_1)\le J_2(\xi_1,\tau_1):=\frac{\lan
\xi_1 \ran^{2|k|-2s}}{\lan\tau_1\ran^{2b}}
\int_{-\infty}^{+\infty}\frac{1} {\lan
\tau_1-\frac{1}{2}\xi_1^2+\xi\xi_1\ran^{2a}} \chi_{\Omega_2}d\xi,
\end{equation}
and
\begin{equation}\label{Estimative-W3-Tilde-a}
\widetilde{W}_3(\xi_2,\tau_2)\le
J_3(\xi_2,\tau_2):=\frac{1}{\lan\tau_2-\tfrac{1}{2}\xi_2^2\ran^{2b}}
\int_{-\infty}^{+\infty}\frac{\lan \xi_1 \ran^{2|k|-2s}} {\lan
\tau_2-\frac{1}{2}\xi_1^2-\frac{1}{2}\xi_2^2 +
\xi_1\xi_2\ran^{2a}} \chi_{\widetilde{\Omega}_3}d\xi_1.
\end{equation}

We begin estimating $J_1$ on $\Omega_1={\mathcal A}_1\cup
{\mathcal A}_2\cup {\mathcal A}_{3,1}$. In region ${\mathcal
A}_1$, using $|\xi_1|\le 1$, $a>0$, $b>1/2$ it easy to see that
\begin{equation}
|J_1|\le C \int_{|\xi_1|\le 1}d\xi_1 \le C.
\end{equation}

In region ${\mathcal A}_2$, by the change of variables $\eta =
\tau +\frac{1}{2}\xi^2+\frac{1}{2}\xi_1^2-\xi\xi_1$ and the
condition $|\xi-\xi_1|\ge \frac{1}{8}|\xi_1|$ we obtain
\begin{equation}\begin{split}
|J_1|&\le \frac{1}{\lan \tau + \frac{1}{2}\xi^2 \ran^{2a}}
\int_{{\mathcal A}_2}\frac{\lan \xi_1 \ran^{2|k|-2s}}
{|\xi_1-\xi|\lan \eta
\ran^{2b}}d\eta\\
&\le \frac{8}{\lan \tau + \frac{1}{2}\xi^2
\ran^{2a}}\int_{{\mathcal A}_2}\frac{\lan \xi_1
\ran^{2|k|-2s}}{|\xi_1|\lan \eta \ran^{2b}}d\eta\\
&\le C,
\end{split}\end{equation}
where we have used that $a>0$, $b>1/2$ and $|k|-s\le 1/2$.

In region ${\mathcal A}_{3,1}$, by (\ref{Algebraic-Inequality}) we
have that $$|\xi_1|^2\leq 24\lan \tau+ \tfrac{1}{2}\xi^2\ran$$ and
consequently using $a>0$ we obtain $$\lan \tau+
\tfrac{1}{2}\xi^2\ran^{-2a}\le C|\xi_1|^{-4a}.$$ Then we use that
$|k|-s \le 1/2<2a$, for $a > 1/4$, combined with Lemma
\ref{pe-sd-lemma2}-(\ref{pe-sd-lemma2-c}) to get
\begin{equation}
|J_1|\le C \int_{\R}\frac{\lan \xi_1 \ran^{2|k|-2s}}
{|\xi_1|^{4a}\lan \tau
+\frac{1}{2}\xi^2+\frac{1}{2}\xi_1^2-\xi\xi_1\ran^{2b}} d\xi_1\le
C.
\end{equation}

Next we estimate $J_2$. First, we making the change
$$\eta= \tau_1-\tfrac{1}{2}\xi_1^2+\xi\xi_1,\;\;\;\;\;\;\;\; d\eta=\xi_1d\xi,$$
and we note that the relations in (\ref{Algebraic-Relation}) and
the restriction in region $\Omega_2$ yield
\begin{equation}\label{J2-A}
\lan \eta \ran \le \lan \tau_1\ran  + |\xi
\xi_1-\tfrac{1}{2}\xi_1^2|\le 4\lan \tau_1\ran.
\end{equation}
Moreover, by (\ref{Algebraic-Inequality}) we have
$$|\xi_1|^2\le 24\lan \tau_1\ran$$
and hence using that $2a+2b-1>0$ we get
\begin{equation}\label{J2-B}
|\xi_1|^{4a+4b-2}\le C\lan \tau_1\ran ^{2a+2b-1}.
\end{equation}
Now using the inequalities (\ref{J2-A}), (\ref{J2-B}) and that $a<
1/2$ we can estimate $J_2$ as follows:

\begin{equation}\begin{split}
|J_2(\xi_1,\tau_1)|&\le \frac{\lan \xi_1\ran^{2|k|-2s}}{\lan
\tau_1\ran^{2b}}\int_{\lan \eta \ran\le 4\lan \tau_1
\ran}\frac{d\eta}{|\xi_1|(1+|\eta|)^{2a}}\\
& \le C\frac{\lan \xi_1 \ran^{2|k|-2s}} {\lan
\tau_1\ran^{2b}|\xi_1|}\lan \tau_1 \ran ^{1-2a}\\
& \le C \frac{\lan\xi_1 \ran^{2|k|-2s}}{\lan\tau_1
\ran^{2a+2b-1}|\xi_1|}\\
& \le C \frac{\lan\xi_1 \ran^{2|k|-2s}}{\lan\xi_1
\ran^{4a+4b-2}|\xi_1|}\\
&\le C,
\end{split}\end{equation}
where the last inequality follows directly from the conditions
$2a+2b-1/2\ge 1/2$ (for $a>0$) and $|k|-s \le 1/2$.

Finally, in region $\tilde{\Omega}_3$ we note that
$$|\xi_1|^{4b}\le C\lan \tau_2-\frac{1}{2}\xi_2^2\ran ^{2b}.$$
Hence, from conditions $a>1/4$,\; $b>1/2$\; and\; $|k|-s \le 1/2$
coupled with Lemma \ref{pe-sd-lemma2}-(\ref{pe-sd-lemma2-c}), we
have that

\begin{equation}\begin{split}
|J_3(\xi_2,\tau_2)|&\le C\int_{\tilde{\Omega}_3}\frac{\lan \xi_1
\ran ^{2|k|-2s}}{|\xi_1|^{4b}\lan \tau_2-\frac{1}{2}\xi_1^2
-\frac{1}{2}\xi_2^2 + \xi_1\xi_2\ran^{2a}}d\xi_1\\
& \le C,
\end{split}\end{equation}
which complete the proof of desired estimate.
\end{proof}

For later use, we note that the following result is a consequence of the proof of the previous proposition:

\begin{corollary}\label{corollary-uv-continuous}It holds $\|uv\|_{X^{0,-1/4+}}\lesssim \|u\|_{X^{0,0}}\|v\|_{H_t^{1/2+}L_x^2}$.
\end{corollary}

\begin{proof}Putting $k=s=0$ in the proof of the proposition~\ref{proposition-uv-continuous}, we see that our task is reduced to show that the following expression is bounded
\begin{equation*}
\widetilde{W}_3:=\sup\limits_{\xi_2,\tau_2}\widetilde{W}_3(\xi_2,\tau_2):= \sup\limits_{\xi_2,\tau_2}\frac{1}{\langle\tau_2-\tfrac{1}{2}\xi_2^2\rangle^{2b_1}}\int_{\mathbb{R}^2}\frac{d\xi_1 d\tau_1}{\langle \tau_1 \rangle^{2b_2} \langle\tau_1-\tau_2+\tfrac{1}{2}(\xi_1-\xi_2)^2\rangle^{2a}}
\end{equation*}
where $a=1/4+$, $b_1=0$ and $b_2=1/2+$. On the other hand, we can use the lemma~\ref{pe-sd-lemma2} in order to obtain that
\begin{equation*}
\widetilde{W}_3\leq \int_{\mathbb{R}^2}\frac{d\xi_1 d\tau_1}{\langle \tau_1 \rangle^{2b_2} \langle\tau_1-\tau_2+\tfrac{1}{2}(\xi_1-\xi_2)^2\rangle^{2a}}\leq \int_{\mathbb{R}}\frac{d\xi_1}{\langle -\tau_2+\tfrac{1}{2}(\xi_1-\xi_2)^2\rangle^{2a}}\lesssim 1.
\end{equation*}
The desired corollary follows.
\end{proof}


\begin{proposition}\label{proposition-|u|2-continuous}
If $\max\{0,s\}< 2k + 1/2$ and $s\leq k+1/2$, then the bilinear
estimate
\begin{equation}\label{Estimativa-Mixta-2-c}
\|u\bar{w}\|_{H^{-a}_tH^s_x}\lesssim
\|u\|_{X^{k,b}}\|w\|_{X^{k,b}}
\end{equation}
holds if $b>1/2$ and $\max\{\frac{1}{4},\max\{0,s\}-2k\}< a <1/2$.
\end{proposition}

\begin{proof}
Analogously to the previous proposition, the estimate
(\ref{Estimativa-Mixta-2-c}) is equivalent to prove
\begin{equation}\label{Desigualdade-Equivalente-2}
 Z_{f,g}(\varphi)
\lesssim \|f\|_{L^2}\|g\|_{L^2}\|\varphi\|_{L^2},
\end{equation}
for all $\varphi \in L^2(\R^2)$, where
\begin{equation}\label{Desigualdade-Equivalente-2-a}
 \small{Z_{f,g}(\varphi)=\int_{\R^4} \frac{\lan \xi
\ran^s\bar{\varphi}(\xi, \tau) f(\xi-\xi_1, \tau
-\tau_1)\bar{g}(-\xi_1, -\tau_1)}{\lan \tau \ran^a\lan \xi-\xi_1
\ran^k \lan \tau -\tau_1 + \frac{1}{2}(\xi-\xi_1)^2\ran^b \lan
\xi_1 \ran^k\lan\tau_1 -\frac{1}{2}\xi_1^2\ran^b}d\xi_1
d\tau_1d\xi d\tau}
\end{equation}

We have the following dispersion relation

\begin{equation}\label{eq-u2-dispersion}
\begin{cases}
\xi=\xi_1+\xi_2,\quad  \tau=\tau_1+\tau_2,\\
\sigma_1=\tau_1-\frac{1}{2}\xi_1^2, \quad \sigma_2=\tau_2+\frac{1}{2}\xi_2^2,\\
\tau -\sigma_1 - \sigma_2=-\tfrac{1}{2}\xi^2 +\xi \xi_1
=\tfrac{1}{2}\xi^2 -\xi \xi_2=\tfrac{1}{2}(\xi_1^2-\xi_2^2).
\end{cases}
\end{equation}

\smallskip
We divide $\R^4$ in the following integration regions:

\smallskip
\textbf{Region $\boldsymbol{A}$:} $\boldsymbol{|\sigma_1|\geq \max
\{|\tau|,\; |\sigma_2|\}}$. We consider two subregions of $A$:

\bigskip
\noindent{\textbf{Subregion $\boldsymbol{A_1}$:}}
$\boldsymbol{|\xi_1|\leq 2|\xi_2|}$. If $k\leq 0$, we have
$\lan\xi\ran^s\lan\xi_1\ran^{-k}\lan\xi_2\ran^{-k}\lesssim
\lan\xi_2\ran^{\max\{0,s\}-2k}\lesssim\langle\xi_2\rangle^{1/2}$ (because
$|\xi|\leq 3|\xi_2|$ and $\max\{0,s\}\leq 2k+1/2$). Hence, we can
estimate
\begin{equation}\label{e.Z-regiao-a1}
 Z\lesssim \int_{\R^4} \frac{\lan \xi_2
\ran^{1/2}\bar{\varphi}(\xi, \tau) f(\xi_2, \tau_2)\bar{g}(-\xi_1,
-\tau_1)}{\lan \tau \ran^a\lan \sigma_2\ran^b \lan
\sigma_1\ran^b}\chi_{A_1} d\xi_1 d\tau_1d\xi d\tau
\end{equation}
if $k\leq 0$. Thus, in the same way as the previous estimate of
(\ref{Estimative-W1}), it suffices to bound the expression:
\begin{equation}\label{e.z1}
\widetilde{Z}_1:=
\sup_{\xi_1,\tau_1}\frac{1}{\lan\sigma_1\ran^{2b}}\int_{\mathbb{R}^2}
\frac{\lan\xi_2\ran\chi_{A_1}}{\lan\tau\ran^{2a}\lan\sigma_2\ran^{2b}}
\end{equation}
\begin{itemize}
\item If $|\xi_2|\le 1$ we using Lemma \ref{pe-sd-lemma2}-(\ref{pe-sd-lemma2-c})
to  get
\begin{equation*}\begin{split}
\widetilde Z_1 &\lesssim \int_{\R^2} \frac{\chi_{A_1}}{\lan \tau
\ran^{2a}\lan \tau -\tau_1 +
\frac{1}{2}(\xi-\xi_1)^2\ran^{2b}}d\xi d\tau\\
&\lesssim \int_{-\infty}^{+\infty}\frac{1}{\lan
-\tau_1+\frac{1}{2}\xi_1^2-\xi\xi_1+\frac{1}{2}\xi^2\ran^{2a}}d\xi\\
&\lesssim 1,
\end{split}
\end{equation*}
since $2a> 1/2.$

\item If $|\xi_2| > 1$ we have that $\lan \xi_2 \ran\lesssim
|\xi_2|$. Next, changing variables  $\tau=\tau_2+\tau_1$ and
$\sigma_2=\tau_2+\tfrac{1}{2}\xi_2^2$, for fixed $\xi_1$ and
$\tau_1$, we have that $d\tau d\sigma_2=|\xi_2|d\tau_2d\xi_2$ and
then we obtain

\begin{equation*}\begin{split}
\widetilde Z_1 &\lesssim \sup_{\xi_1, \tau_1}\frac{1}{\lan \sigma_1
\ran ^{2b}} \int_{\R^2} \frac{|\xi_2|\chi_{A_1}}{\lan \tau
\ran^{2a}\lan \sigma_2\ran^{2b}}d\xi_2 d\tau_2\\
&=\sup_{\xi_1, \tau_1}\frac{1}{\lan \sigma_1 \ran ^{2b}}\int_{\R^2}
\frac{\chi_{A_1}}{\lan \tau
\ran^{2a}\lan \sigma_2\ran^{2b}}d\tau d\sigma_2 \\
& \lesssim \sup_{\xi_1, \tau_1}\frac{1}{\lan \sigma_1 \ran ^{2b}}
\int_{0}^{|\sigma_1|}\lan \tau \ran^{-2a}d\tau\int_{0}^{|\sigma_1|}
\lan \sigma_2\ran^{-2b}d\sigma_2\\
& \lesssim \sup_{\sigma_1}\lan \sigma_1\ran^{-2b}\lan
\sigma_1\ran^{1-2a}\lan \sigma_1\ran^{1-2b}\\
&=\sup_{\sigma_1}\lan \sigma_1 \ran^{2-2a-4b}\\
&\lesssim 1,
\end{split}
\end{equation*}
since $0<a$ and $b>1/2$ implies $2-2a-4b\leq 0$.
\end{itemize}

Therefore, we showed (\ref{Desigualdade-Equivalente-2}) in the
subregion $A_1$ whenever $k\leq 0$. On the other, if $k\geq 0$, we
have $\lan\xi\ran^s\lan\xi_1\ran^{-k}\lan\xi_2\ran^{-k}\lesssim
\lan\xi\ran^{s-k}\lesssim\lan\xi\ran^{1/2}$ (because $|\xi|\leq
3|\xi_2|$ and $s-k\leq 1/2$). So, we get
\begin{equation}
 Z\lesssim \int_{\R^4} \frac{\lan \xi
\ran^{1/2}\bar{\varphi}(\xi, \tau) f(\xi_2, \tau_2)\bar{g}(-\xi_1,
-\tau_1)}{\lan \tau \ran^a\lan \sigma_2\ran^b \lan
\sigma_1\ran^b}\chi_{A_1} d\xi_1 d\tau_1d\xi d\tau
\end{equation}
if $k\geq 0$. Thus, applying the lemma~\ref{l.ginibre}, we also
obtain (\ref{Desigualdade-Equivalente-2}) if $k\geq 0$.

This completes the analysis of $Z$ in the subregion
$\boldsymbol{A_1}$.

\bigskip
\noindent{\textbf{Subregion $\boldsymbol{A_2}$:}}
$\boldsymbol{|\xi_1|\ge 2|\xi_2|}$. Here, the dispersion relation
(\ref{eq-u2-dispersion}) yields that
$$\tfrac{3}{4}\xi_1^2\le |\xi_1^2-\xi_2^2|=2|\tau -\sigma_1-\sigma_2|\le 6|\sigma_1|
\Longrightarrow \xi_1^2\le 8|\sigma_1|.$$ Hence,
\begin{equation}\label{Estimate-Sigma1-Region-A2}
\frac{1}{\lan \sigma_1 \ran}\lesssim \frac{1}{\lan \xi_1 \ran^2}.
\end{equation}
If $k\leq 0$, it follows $\lan\xi\ran^s \lan\xi_1\ran^{-k}
\lan\xi_2\ran^{-k} \lesssim \lan\xi_1\ran^{1/2}$ (because
$\max\{0,s\}\leq 2k+1/2$ and $|\xi|\leq 3|\xi_1|/2$), so that
\begin{equation}
Z\lesssim \int_{\R^4} \frac{\lan \xi_1 \ran^{1/2}\bar{\varphi}(\xi,
\tau) f(\xi_2, \tau_2)\bar{g}(-\xi_1, -\tau_1)}{\lan \tau \ran^a\lan
\sigma_2\ran^b \lan \sigma_1\ran^b}\chi_{A_2} d\xi_1 d\tau_1d\xi
d\tau
\end{equation}
if $k\leq 0$. Thus, similarly to (\ref{e.Z-regiao-a1}), our task is
to estimate
\begin{equation}
\widetilde{Z}_2:=\sup_{\xi_1,\tau_1} \frac{1}{\lan\sigma_1\ran^{2b}}
\int_{\mathbb{R}^2}\frac{\lan\xi_1\ran}{\lan\tau\ran^{2a}\lan\sigma_2\ran^{2b}}\chi_{A_2}
d\xi_2 d\tau_2
\end{equation}
Using (\ref{Estimate-Sigma1-Region-A2}), lemma
\ref{pe-sd-lemma2}-(\ref{pe-sd-lemma2-a}) and lemma
\ref{pe-sd-lemma2}-(\ref{pe-sd-lemma2-c}) we obtain
\begin{equation}\begin{split}\label{Desigualdade-Equivalente-2-g}
\widetilde Z_2& \lesssim \sup_{\xi_1,\tau_1}{\lan \xi_1 \ran
^{1-4b}} \int_{\R^2} \frac{\chi_{A}}{\lan \tau \ran^{2a}\lan \tau
-\tau_1
+ \frac{1}{2}(\xi-\xi_1)^2\ran^{2b}}d\xi d\tau\\
& \lesssim \sup_{\xi_1, \tau_1}{\lan \xi_1 \ran^{1-4b}}
\int_{-\infty}^{+\infty}\frac{1}{\lan
-\tau_1+\frac{1}{2}\xi_1^2-\xi\xi_1+\frac{1}{2}\xi^2\ran^{2a}}d\xi\\
& \lesssim 1,
\end{split}\end{equation}
where in the last inequality we have used that since $1/4< a$ and
$b>1/4$.

If $k\geq 0$, we have $\lan\xi\ran^s \lan\xi_1\ran^{-k}
\lan\xi_2\ran^{-k} \lesssim \lan\xi\ran^{s-k}\lesssim
\lan\xi\ran^{1/2}$ since $s-k\leq 1/2$ and $|\xi|\leq 3|\xi_1|/2$.
So, we get
\begin{equation}
Z\lesssim\int_{\R^4} \frac{\lan \xi \ran^{1/2}\bar{\varphi}(\xi,
\tau) f(\xi_2, \tau_2)\bar{g}(-\xi_1, -\tau_1)}{\lan \tau \ran^a\lan
\sigma_2\ran^b \lan \sigma_1\ran^b}\chi_{A_2} d\xi_1 d\tau_1d\xi
d\tau\lesssim \|f\|_{L^2}\|g\|_{L^2}\|\varphi\|_{L^2}
\end{equation}
by lemma~\ref{l.ginibre}. This completes the analysis of the $Z$ in
the subregion $\boldsymbol{A_2}$.

\bigskip

Clearly $\boldsymbol{A}=\boldsymbol{A_1}\cup\boldsymbol{A_2}$, so
that the estimate (\ref{Desigualdade-Equivalente-2}) holds true in
the region $\boldsymbol{A}$.

\bigskip \textbf{Region $\boldsymbol{B}$:} $\boldsymbol{|\sigma_2|\ge
\max \{|\tau|,\; |\sigma_1|\}}$. The computations for this region
can be obtained from the previous ones (in region $A$) since all the
involved expressions are symmetric under the exchange of the indices
1 and 2.

\bigskip \textbf{Region $\boldsymbol{C}$:} $\boldsymbol{|\tau|\ge
\max \{|\sigma_1|,\; |\sigma_2|\}}$. Here, we analyze several cases
for the frequencies $\xi$ and $\xi_1$.

We begin with the high frequencies for $\xi$, that is:

\smallskip
\noindent{\textbf{Subregion $\boldsymbol{C_1}$:}} $\boldsymbol{|\xi|
\ge 1}$. We separate this region into two smaller subregions.

\smallskip
\noindent{\textbf{Subregion $\boldsymbol{C_{1,1}}$:}}
$\boldsymbol{\left|\xi_1-\frac{1}{2}\xi  \right| \le 1}$. Here we
have that
$$|\xi_1|\le \left| \xi_1 - \tfrac{1}{2}\xi\right| + \left|\tfrac{1}{2}\xi
\right|\Longrightarrow \lan \xi_1 \ran \lesssim \lan \xi \ran$$
and
$$|\xi_2|=\left |\tfrac{1}{2}\xi + \tfrac{1}{2}\xi -\xi_1\right|\le
 \left| \xi_1 - \tfrac{1}{2}\xi\right| + \left|\tfrac{1}{2}\xi
\right|\Longrightarrow \lan \xi_2 \ran \lesssim \lan \xi \ran.$$

In particular, we get $\lan\xi\ran^s \lan\xi_1\ran^{-k}
\lan\xi_2\ran^{-k}\lesssim \lan\xi\ran^{1/2}$ (because
$\max\{0,s\}\leq 2k+1/2$ and $s-k\leq 1/2$). This allows us to
conclude that
\begin{equation}
Z\lesssim\int_{\R^4} \frac{\lan \xi \ran^{1/2}\bar{\varphi}(\xi,
\tau) f(\xi_2, \tau_2)\bar{g}(-\xi_1, -\tau_1)}{\lan \tau \ran^a\lan
\sigma_2\ran^b \lan \sigma_1\ran^b}\chi_{C_{1,1}} d\xi_1 d\tau_1d\xi
d\tau\lesssim \|f\|_{L^2}\|g\|_{L^2}\|\varphi\|_{L^2}
\end{equation}
by lemma~\ref{l.ginibre}, which is the desired estimate
(\ref{Desigualdade-Equivalente-2}) in the subregion
$\boldsymbol{C_{1,1}}$.

\smallskip
\noindent{\textbf{Subregion $\boldsymbol{C_{1,2}}$:}}
$\boldsymbol{\left|\xi_1-\frac{1}{2}\xi  \right| \ge 1}$. Firstly,
we note that if $\min\{|\xi_1|, |\xi_2|\}\le 1$, it follows that
$\max\{\lan\xi_1\ran, \lan\xi_2\ran\}\lesssim \lan\xi\ran$ and the
same analysis of the subregion $\boldsymbol{C_{1,1}}$ can be
repeated here. Thus, we can assume that $|\xi_1|\ge 1$ and
$|\xi_2|\ge 1$. Note that

\begin{equation}\label{Bilinear-U2-RC1-2-a}
\begin{split}
|\xi_1 \xi_2|=|\xi_1
(\xi-\xi_1)|&=\left|\left((\xi_1-\tfrac{1}{2}\xi)+\tfrac{1}{2}\xi\right)
\left(\tfrac{1}{2}\xi-(\xi_1-\tfrac{1}{2}\xi)\right)\right|\\
& \leq \left|\xi_1-\tfrac{1}{2}\xi\right|^2+\tfrac{1}{4}|\xi|^2.
\end{split}
\end{equation}
Also, from (\ref{eq-u2-dispersion}) and the conditions $|\xi|\ge 1$
and $|\xi_1-\tfrac{1}{2}\xi|\ge1$, it follows that
\begin{equation}\label{Bilinear-U2-RC1-2-b}
\max\{|\xi|, |\xi_1-\tfrac{1}{2}\xi|\} \le
\left|\xi(\xi_1-\tfrac{1}{2}\xi)\right|\le 3\lan \tau \ran.
\end{equation}
If $s\leq 0$, $k\leq 0$, we obtain $\lan\xi\ran^s \lan\xi_1\ran^{-k}
\lan\xi_2\ran^{-k}\lesssim \lan\xi_1\ran^{\tfrac{\max\{0,s\}}{2}-k}
\lan\xi_2\ran^{\tfrac{\max\{0,s\}}{2}-k}$; if $s\leq 0$, $k\geq 0$, we get
$\lan\xi\ran^s \lan\xi_1\ran^{-k} \lan\xi_2\ran^{-k}\lesssim 1$; in
the remaining cases (i.e., either $s\geq 0$, $k\leq 0$ or $s\geq 0$,
$k\geq 0$), we have two possibilities, namely $|\xi_1|\sim |\xi_2|$
or $|\xi_1|\nsim |\xi_2|$; when the first case occurs, it follows
that $\lan\xi\ran^s\lan\xi_1\ran^{-k}\lan\xi_2\ran^{-k}\lesssim
\lan\xi_1\ran^{s-2k}\lesssim\lan\xi_1\ran^{\tfrac{\max\{0,s\}}{2}-k}\lan\xi_2\ran^{\tfrac{\max\{0,s\}}{2}-k}$
and, in the second case, we conclude that
$\lan\tau\ran\gtrsim\lan\xi\ran^2$, which implies
$\lan\xi\ran^s\lan\xi_1\ran^{-k}\lan\xi_2\ran^{-k}\lesssim
\lan\xi\ran^{s-2k}\lesssim
\lan\xi\ran^{1/2}\lesssim\lan\tau\ran^{1/4}$.

In resume, we always get that, in any case,
either $\lan\xi\ran^s\lan\xi_1\ran^{-k}\lan\xi_2\ran^{-k}\lesssim
\lan\tau\ran^{1/4}$ or $\lan\xi\ran^s\lan\xi_1\ran^{-k}\lan\xi_2\ran^{-k}\lesssim \lan\xi_1\ran^{\tfrac{\max\{0,s\}}{2}-k}\lan\xi_2\ran^{\tfrac{\max\{0,s\}}{2}-k}$. When the first possibility occurs, using Cauchy-Schwarz, we can reduce the estimate (\ref{Desigualdade-Equivalente-2}) to bound the
expression:

\begin{equation}
\widetilde{Z}:=\sup_{\xi, \tau}\frac{1}{\lan \tau \ran ^{2a}}
\int_{\R^2} \frac{\lan \tau \ran^{1/2}\chi_{C_{1,2}}}{\lan \sigma_1
\ran^{2b}\lan \sigma_2\ran^{2b}}d\xi_2 d\tau_2.
\end{equation}
But, this can be done as follows:
\begin{equation}\begin{split}\label{Bilinear-U2-RC1-2-c}
\widetilde{Z} &\lesssim \sup_{\xi, \tau}\frac{\lan \tau
\ran^{1/2-2a}}{ |\xi|}\int_{\R^2} \frac{ \chi_{C_{1,2}}}{\lan
\sigma_1
\ran^{2b}\lan \sigma_2\ran^{2b}}d\sigma_1 d\sigma_2 \\
& \lesssim \sup_{\xi, \tau}\frac{\lan \tau\ran^{1/2-2a}}{|\xi|}\lan \tau \ran^{2-4b}\\
&\lesssim 1,
\end{split}
\end{equation}
since $|\xi|\geq 1$, $b>1/2$ and $a>1/4$. When the second possibility happens, we decompose the frequencies $\xi_j$ and the modulations $\sigma_j$ into dyadic blocks $\lan\xi_j\ran\sim N_j$ and $\lan\sigma_j\ran\sim L_j$ (here $\xi_0:=\xi$, $\sigma_0:= \tau$ and $j=0,1,2$). Hence, it suffices to estimate (\ref{Desigualdade-Equivalente-2}) restricted to each dyadic block with the gain of extra terms $N_j^{0-}$ and $L_j^{0-}$. To simplify, we put $N_{\max}:=\max\{N_0,N_1,N_2\}$ and $L_{\max}:= \max\{L_0,L_1,L_2\}$. So, we have

\begin{equation*}
Z\lesssim\int_{\R^4} \frac{\lan\xi_1\ran^{\tfrac{\max\{0,s\}}{2}-k}\lan\xi_2\ran^{\tfrac{\max\{0,s\}}{2}-k}\bar{\varphi}(\xi,
\tau) f(\xi_2, \tau_2)\bar{g}(-\xi_1, -\tau_1)}{\lan \tau \ran^a\lan
\sigma_2\ran^b \lan \sigma_1\ran^b}\chi_{C_{1,2}} d\xi_1 d\tau_1d\xi
d\tau
\end{equation*}
Using (\ref{Bilinear-U2-RC1-2-a}) and (\ref{eq-u2-dispersion}), we get $\lan\xi_1\ran\lan\xi_2\ran\lesssim \lan\tau\ran^2$. Since $a>\max\{0,s\}-2k$, we get
\begin{equation*}
Z\lesssim\int_{\R^4} \frac{\bar{\varphi}(\xi,
\tau) f(\xi_2, \tau_2)\bar{g}(-\xi_1, -\tau_1)}{L_0^{0+}\lan
\sigma_2\ran^b \lan \sigma_1\ran^b}\chi_{C_{1,2}}
\end{equation*}
Applying Cauchy-Schwarz and the lemma~\ref{l.ginibre}, it suffices to bound the expression:
\begin{equation*}
\sup_{\xi, \tau}\frac{1}{L_0^{0+}}
\int_{\R^2} \frac{\chi_{C_{1,2}}}{\lan \sigma_1
\ran^{2b}\lan \sigma_2\ran^{2b}}d\xi_2 d\tau_2
\end{equation*}
Recall that (\ref{Bilinear-U2-RC1-2-a}), (\ref{Bilinear-U2-RC1-2-b}) and (\ref{eq-u2-dispersion}) implies $N_{\max}\lesssim L_0$. Also, $L_0=L_{\max}$ in the region $\boldsymbol{C}$. In particular,
\begin{equation*}
\begin{split}
\sup_{\xi, \tau}\frac{1}{L_0^{0+}}
\int_{\R^2} \frac{\chi_{C_{1,2}}}{\lan \sigma_1
\ran^{2b}\lan \sigma_2\ran^{2b}}d\xi_2 d\tau_2&\lesssim
\sup_{\xi, \tau} \frac{N_{\max}^{0-} L_{\max}^{0-}}{|\xi|} \int_{\R^2} \frac{\chi_{C_{1,2}}}{\lan \sigma_1
\ran^{2b}\lan \sigma_2\ran^{2b}}d\sigma_1 d\sigma_2 \\
&\lesssim\sup_{\xi, \tau} \frac{N_{\max}^{0-} L_{\max}^{0-}}{|\xi|} \lan\tau\ran^{2-4b} \\
&\lesssim N_{\max}^{0-} L_{\max}^{0-},
\end{split}
\end{equation*}
because $b>1/2$ and $|\xi|\geq 1$. This completes our analysis of the region $\boldsymbol{C_{1,2}}$.

We conclude with the small frequencies for $\xi$, that is:

\noindent{\textbf{Subregion $\boldsymbol{C_2}$:}} $\boldsymbol{|\xi|
\le 1}$. The hypothesis $|\tau|\ge \max \{|\sigma_1|,\;
|\sigma_2|\}$ is not crucial in this case; hence we divide into two smaller
subregions:

\smallskip
\noindent{\textbf{Subregion $\boldsymbol{C_{2,1}}$:}}
$\boldsymbol{\left|\xi_1 \right|\le 2}$. Here, it is easy to see
that $\lan \xi_1\ran\lesssim1$, $\lan \xi_2 \ran \lesssim 1$. In
particular, by Cauchy-Schwarz, our task is to estimate
\begin{equation*}
\widetilde{Z}=\sup_{\xi_1, \sigma_1}
\frac{1}{\lan\sigma_1\ran^{2b}} \int_{\R^2} \frac{1}{\lan \tau \ran^{2a}\lan
\sigma_2\ran^{2b}
}d\xi_2 d\tau_2.
\end{equation*}
Then, using lemma \ref{pe-sd-lemma2}-(\ref{pe-sd-lemma2-a}) and lemma \ref{pe-sd-lemma2}-(\ref{pe-sd-lemma2-c}), we get
\begin{equation}\begin{split}\label{Bilinear-U2-RC2-1}
\widetilde Z&=\sup_{\xi_1, \sigma_1}\frac{1}{\lan \sigma_1 \ran
^{2b}} \int_{\R^2} \frac{1}{\lan \tau \ran^{2a}\lan
\sigma_2\ran^{2b}
}d\xi_2 d\tau_2\\
&\lesssim \int_{\R^2} \frac{1}{\lan \tau \ran^{2a}\lan \tau
-\tau_1 +
\frac{1}{2}(\xi-\xi_1)^2\ran^{2b}}d\xi d\tau\\
&\lesssim \int_{-\infty}^{+\infty}\frac{1}{\lan
-\tau_1+\frac{1}{2}\xi_1^2-\xi\xi_1+\frac{1}{2}\xi^2\ran^{2a}}d\xi\\
&\lesssim 1,
\end{split}\end{equation}
since $a>1/4.$

\smallskip
\noindent{\textbf{Subregion $\boldsymbol{C_{2,2}}$:}}
$\boldsymbol{\left |\xi_1 \right|\ge 2}$.
Redoing the analysis of the bounds for the term $\lan\xi\ran^s\lan\xi_1\ran^{-k}\lan\xi_2\ran^{-k}$ in the four cases $s\leq 0$, $k\leq 0$, $\dots$, $s\geq 0$, $k\geq 0$, we see that
\begin{equation*}
\lan\xi\ran^s\lan\xi_1\ran^{-k}\lan\xi_2\ran^{-k}\lesssim \lan\xi_1\ran^{1/2}.
\end{equation*}
Similarly to the previous estimates for subregion $C_{1,2}$, we decompose the frequencies $\lan\xi_j\ran\sim N_j$, $j=0,1,2$, into dyadic blocks so that our task is to bound (\ref{Desigualdade-Equivalente-2}) restricted to each dyadic block with the gain of extra terms $N_j^{0-}$. We have
\begin{equation*}
Z\lesssim N_1^{1/2}\int_{\R^4}\frac{\bar{\varphi}(\xi,
\tau) f(\xi_2, \tau_2)\bar{g}(-\xi_1, -\tau_1)}{\lan\tau\ran^a\lan
\sigma_2\ran^b \lan \sigma_1\ran^b}\chi_{C_{2,2}} d\xi_1 d\tau_1d\xi
d\tau
\end{equation*}
Applying Cauchy-Schwarz, it suffices to prove that:
\begin{equation}\label{e.|u|2-dyadic2}
N_1^{1/2}\sup_{\xi_1, \sigma_1}\frac{1}{\lan \sigma_1 \ran
^{2b}}\int_{\R^2}\frac{\chi_{C_{2,2}}}{\lan \tau \ran^{2a}\lan \tau
-\tau_1 + \frac{1}{2}(\xi-\xi_1)^2\ran^{2b}}d\xi d\tau\lesssim N_{\max}^{0-}.
\end{equation}
This can be accomplished as follows. Firstly, notice that
\begin{equation}\begin{split}\label{Bilinear-U2-RC2-2}
&\sup_{\xi_1, \sigma_1}\frac{1}{\lan \sigma_1 \ran
^{2b}}\int_{\R^2}\frac{\chi_{C_{2,2}}}{\lan \tau \ran^{2a}\lan \tau
-\tau_1 + \frac{1}{2}(\xi-\xi_1)^2\ran^{2b}}d\xi d\tau\\
&\lesssim \sup_{\xi_1, \sigma_1}\frac{1}{\lan
\sigma_1 \ran^{2b}}\int_{|\xi| \le 1}\frac{d\xi}{\lan
-\tau_1+\frac{1}{2}\xi_1^2-\xi\xi_1+\frac{1}{2}\xi^2\ran^{2a}}.
\end{split}\end{equation}
Now, by  changing variables
$$\eta = -\tau_1+\frac{1}{2}\xi_1^2-\xi\xi_1+\frac{1}{2}\xi^2,\quad d\eta =(\xi-\xi_1)d\xi,$$
we get $|\eta|\le \lan \sigma_1\ran + |\xi_1| + \frac{1}{2} \le  \lan \sigma_1\ran + 2|\xi_1|$ and
we obtain the following bound of (\ref{Bilinear-U2-RC2-2}):

\begin{equation}\begin{split}\label{Bilinear-U2-RC2-2-a}
&\sup_{\xi_1, \sigma_1}
\frac{1}{\lan \sigma_1 \ran^{2b}}\int_{|\eta| \leq \lan \sigma_1\ran + 2|\xi_1|}\frac{d\eta}{(1+|\eta|)^{2a}|\xi_1-\xi|}\\
&\lesssim \sup_{\xi_1, \sigma_1}\frac{1}{\lan \sigma_1 \ran
^{2b}|\xi_1|}\int_{|\eta| \le \lan \sigma_1\ran + 2|\xi_1|}\frac{d\eta}{(1+|\eta|)^{2a}}\\
&\lesssim \sup_{\xi_1, \sigma_1}\frac{1}{\lan \sigma_1 \ran
^{2b}|\xi_1|}\left(\lan \sigma_1 \ran
^{1-2a} + |\xi_1|^{1-2a}\right )\\
& \lesssim \frac{1}{|\xi_1|^{2a}},
\end{split}\end{equation}
since $2b> 1$ and $a>0$. Putting this estimate into the expression (\ref{e.|u|2-dyadic2}), because $a>1/4$ and $N_1\sim N_{\max}$ in the subregion $\boldsymbol{C_{2,2}}$, we conclude
\begin{equation}
\begin{split}
&N_1^{1/2}\sup_{\xi_1, \sigma_1}\frac{1}{\lan \sigma_1 \ran
^{2b}}\int_{\R^2}\frac{\chi_{C_{2,2}}}{\lan \tau \ran^{2a}\lan \tau
-\tau_1 + \frac{1}{2}(\xi-\xi_1)^2\ran^{2b}}d\xi d\tau\\
&\lesssim N_1^{1/2}\cdot \frac{1}{N_1^{2a}} \\
&\lesssim N_{\max}^{0-}.
\end{split}
\end{equation}

\smallskip
Collecting all the estimates above in all regions we have that the inequality
(\ref{Desigualdade-Equivalente-2}) holds provided the conditions
in proposition~\ref{proposition-|u|2-continuous} are valid.
\end{proof}

\begin{corollary}\label{corollary-|u2|-continuous}It holds $\|u\bar{w}\|_{H_t^{-1/4+}L_x^2}\lesssim \|u\|_{X^{0,1/2+}} \|w\|_{X^{0,0}}$.
\end{corollary}

\begin{proof}From the proof of the previous proposition with $k=s=0$, we know that it suffices to show that
\begin{equation*}
\widetilde{Z}_1:=\sup\limits_{\xi_1,\tau_1}\frac{1}{\langle\sigma_1\rangle^{2b_1}}\int\limits_{\mathbb{R}^2} \frac{d\xi d\tau}{\langle\tau\rangle^{2a} \langle\sigma_2\rangle^{2b_2}}\lesssim 1
\end{equation*}
where $a=1/4+$, $b_1=0$ and $b_2=1/2+$. However, this is a simple application of the lemma~\ref{pe-sd-lemma2}:
\begin{equation*}
\widetilde{Z}_1\lesssim\int\limits_{\mathbb{R}^2} \frac{d\xi d\tau}{\langle\tau\rangle^{2a} \langle\sigma_2\rangle^{2b_2}}\lesssim \int\limits_{\mathbb{R}} \frac{d\xi}{\langle-\tau_1+\frac{1}{2}(\xi-\xi_1)^2\rangle^{2a}}\lesssim 1.
\end{equation*}
This ends the argument.
\end{proof}

\begin{remark}As pointed out in the introduction, once the bilinear estimates in propositions \ref{proposition-uv-continuous} and \ref{proposition-|u|2-continuous} are established, it is a standard matter to conclude the local well-posedness statement of theorem~\ref{local-theorem-continuous}. We refer the reader to the works \cite{KPV}, \cite{Bekiranov1} and \cite{Ginibre} for further details.
\end{remark}

\subsection{Counter-Examples I: the continuous case}\label{Examples-Continuous} We finish this
section exhibiting several counter-examples showing that the
bilinear estimates proved above are sharp, that is, the
conditions imposed on the indices $k$ and $s$ in the propositions
\ref{proposition-uv-continuous} and \ref{proposition-|u|2-continuous} are necessary.

\begin{proposition}\label{Example-Continuous-uv}For any $b_1, b_2 \in \R$,
the estimate $\|uv\|_{X^{k,-1/2}}\lesssim \|u\|_{X^{k,b_1}} \|v\|_{H_t^{b_2} H_x^s}$ holds only if $|k|\leq s+1/2$.
\end{proposition}

\begin{proof}
Take $N\in \Z^+$ a large integer and define
\begin{equation*}
\begin{split}
&A_1= \left \{(\zeta,\eta)\in \R^2;\; 0 \le \zeta \le
1/N \;\; \text{and}\;\; |\eta + \tfrac{1}{2}\zeta^2|\le 1
\right \},\\
& B_1= \left \{(\zeta,\eta)\in \R^2;\; N \le \zeta
\le N+\tfrac{1}{N}\;\; \text{and}\;\; |\eta|\le 1\right \},\\
& A_2= \left \{ (\zeta,\eta)\in \R^2;\; N\le \zeta \le
N+\tfrac{1}{N}\;\;
\text{and}\;\; |\eta+ \tfrac{1}{2}\zeta^2| \le 1 \right \},\\
& B_2= \left \{(\zeta,\eta)\in \R^2;\; -N \le \zeta \le -N+
\tfrac{1}{N} \;\; \text{and}\;\; |\eta|\le 1
\right \}.
\end{split}
\end{equation*}
Put $\widehat{f_j}(\zeta,\eta):=\chi_{A_j}$ and $\widehat{g_j}(\zeta,\eta):=\chi_{B_j}$. A straightforward computation gives that
\begin{equation*}
\begin{split}
& \|f_1g_1\|_{X^{k,-1/2}}\sim \left(\frac{1}{N}\left(\frac{N^k}{N}\right)^2\right)^{1/2} \sim N^{k-\tfrac{3}{2}},\\
& \|f_1\|_{X^{k,b_1}}\sim N^{-1/2} \quad \text{and}\quad \|g_1\|_{H_t^{b_2} H_x^s}\sim N^{s-\tfrac{1}{2}}.
\end{split}
\end{equation*}
So, $\|f_1g_1\|_{X^{k,-1/2}}\lesssim \|f_1\|_{X^{k,b_1}} \|g_1\|_{H_t^{b_2} H_x^s}$ implies that $k\leq s+\tfrac{1}{2}$.
Analogously, another simple computation shows that
\begin{equation*}
\|f_2g_2\|_{X^{k,-1/2}}\sim N^{-\tfrac{3}{2}},\quad
\|f_2\|_{X^{k,b_1}}\sim N^{k-\tfrac{1}{2}}\quad  \text{and}\quad \|g_2\|_{H_t^{b_2} H_x^s}\sim N^{s-\tfrac{1}{2}}.
\end{equation*}
Thus, $\|f_2g_2\|_{X^{k,-1/2}}\lesssim \|f_2\|_{X^{k,b_1}} \|g_2\|_{H_t^{b_2} H_x^s}$ implies that $-k\leq s+1/2$. This completes the proof of the proposition.
\end{proof}

\begin{proposition}\label{Example-Continuous-|u|2}For any $b_1, b_2 \in \R$,
the estimate $\|u\bar{w}\|_{H_t^{-1/2} H_x^s}\lesssim \|u\|_{X^{k,b_1}} \|w\|_{X^{k,b_2}}$ holds only if $s\leq k+1/2$ and $\max\{0,s\}\leq 2k+1/2$.
\end{proposition}

\begin{proof}
Take $N\in \Z^+$ a large integer and define
\begin{equation*}
\begin{split}
& A_1= \left \{(\zeta,\eta)\in \R^2;\; 0 \le \zeta \le
1/N \;\; \text{and}\;\; |\eta + \tfrac{1}{2}\zeta^2|\le 1
\right \},\\
& B_1= \left \{(\zeta,\eta)\in \R^2;\; N \le \zeta
\le N+\tfrac{1}{N}\;\; \text{and}\;\; |\eta
+ \tfrac{1}{2}\zeta^2|\le 1\right \},\\
& A_2= \left \{ (\zeta,\eta)\in \R^2;\; N\le \zeta \le
N+\tfrac{1}{N}\;\;
\text{and}\;\; |\eta+ \tfrac{1}{2}\zeta^2| \le 1 \right \},\\
& B_2= \left \{(\zeta,\eta)\in \R^2;\; -N \le \zeta \le -N+
\tfrac{1}{N} \;\; \text{and}\;\; |\eta + \tfrac{1}{2}\zeta^2|\le 1
\right \},\\
&B_3= \left \{(\zeta,\eta)\in \R^2;\; N \le \zeta
\le N+\tfrac{1}{N}\;\; \text{and}\;\; |\eta
+ \tfrac{1}{2}\zeta^2|\le 1\right \}.
\end{split}
\end{equation*}
Put $\widehat{f_j}(\zeta,\eta):=\chi_{\textsf{A}_j}$ and $\widehat{g_j}(\zeta,\eta):=\chi_{\textsf{B}_j}$ ($j=1,2$). A simple calculation shows that
\begin{equation*}
\begin{split}
& \|f_1\bar{g_1}\|_{H_t^{-1/2}H_x^s}\sim \left(\frac{1}{N}\left(\frac{N^s}{N}\right)^2\right)^{1/2} \sim N^{s-\tfrac{3}{2}},\\
& \|f_1\|_{X^{k,b_1}}\sim N^{-1/2}\quad  \text{and}\quad
\|g_1\|_{X^{k,b_2}}\sim N^{k-\tfrac{1}{2}}.
\end{split}
\end{equation*}
Hence, $\|f_1\bar{g_1}\|_{H_t^{-1/2} H_x^s}\lesssim \|f_1\|_{X^{k,b_1}} \|g_1\|_{X^{k,b_2}}$ implies that $s\leq k+\tfrac{1}{2}$. From the similar way,
we have that
\begin{equation*}
\begin{split}
& \|f_2\bar{g_2}\|_{H_t^{-1/2}H_x^s}\sim \left(\frac{1}{N}\left(\frac{N^s}{N}\right)^2\right)^{1/2} \sim N^{s-\tfrac{3}{2}},\\
& \|f_2\bar{g_3}\|_{H_t^{-1/2}H_x^s}\sim \left(\frac{1}{N}\left(\frac{1}{N}\right)^2\right)^{1/2} \sim N^{-\tfrac{3}{2}},\;\; \text{and}\;\;\\
& \|f_2\|_{X^{k,b_1}}\sim \|g_2\|_{X^{k,b_2}}\sim \|g_3\|_{X^{k,b_2}}\sim N^{k-\tfrac{1}{2}}.
\end{split}
\end{equation*}
Thus, $\|f_2\bar{g_2}\|_{H_t^{-1/2} H_x^s}\lesssim \|f_2\|_{X^{k,b_1}} \|g_2\|_{X^{k,b_2}}$\;
and\; $\|f_2\bar{g_3}\|_{H_t^{-1/2} H_x^s}\lesssim \|f_2\|_{X^{k,b_1}} \|g_3\|_{X^{k,b_2}}$
imply, respectively,  that $s\leq 2k+ \tfrac{1}{2}$ and $0\leq 2k+\tfrac{1}{2}$. Therefore, $\max\{0,s\}\leq 2k+\tfrac{1}{2}$.
\end{proof}

\section{\textbf{Bilinear Estimates  for the Coupling Terms in the Periodic Case}}
\label{section-bilinear-periodic} Here, we show sharp
bilinear estimates for the coupling terms in the periodic setting.

\subsection{Proof of the bilinear estimates II: the periodic case}
\begin{proposition}\label{proposition-uv-periodic}
The bilinear estimate
\begin{equation}\label{Estimativa-uv}
\|uv\|_{X^{k,-\frac{1}{2}+}}\lesssim
\|u\|_{X^{k,\frac{1}{2}-}}\|v\|_{H^{\frac{1}{2}-}_t
H_x^{s}}
\end{equation}
holds if $0\leq s\leq 2k$ and \small{$|k-s|<1$}.
\end{proposition}
\begin{proof}
Fix $s\geq 0$ and $k<s+1$. Taking $a=b=c=1/2-$, our task is to show the bilinear estimate
\begin{equation*}
\|uv\|_{X^{k,-a}}\lesssim \|u\|_{X^{k,b}} \|v\|_{H_t^cH_x^s}
\end{equation*}
Defining $f(n,\tau):=\lan\tau+n^2\ran^b\lan n\ran^k \widehat{u}(n,\tau)$ and $g(n,\tau):= \lan\tau\ran^c \lan n\ran^s
\widehat{v}(n,\tau)$, it suffices to prove that
\begin{equation}\label{e.uv-periodic-1}
Z\lesssim \|f\|_{L^2_{n,\tau}} \|g\|_{L^2_{n,\tau}} \|\varphi\|_{L^2_{n,\tau}},
\end{equation}
where
\begin{equation}\label{e.uv-periodic-2}
W:= \sum\limits_{n\in\mathbb{Z}} \int d\tau \sum\limits_{n=n_1+n_2} \int_{\tau=\tau_1+\tau_2}
\frac{\lan\tau+n^2\ran^{-a} \lan n\ran^k f(n_1,\tau_1) g(n_2,\tau_2) \varphi(n,\tau)}{\lan\tau_1+n_1^2\ran^b
\lan\tau_2\ran^c \lan n_1\ran^k \lan n_2\ran^s}.
\end{equation}
Dividing $\mathbb{Z}^2\times \mathbb{R}^2$ into three regions, namely $\mathbb{Z}^2\times \mathbb{R}^2 = R_0 \cup R_1
\cup R_2$, integrating first over $n_1,\tau_1$ in the region $R_0$, $n,\tau$ in the region $R_1$, $n_2,\tau_2$ in the
region $R_2$ and using Cauchy-Schwarz, we easily see that it remains only to uniformly bound the following three
expressions:

\begin{equation}\label{e.w1}
W_1:=\sup\limits_{n,\tau} \frac{\lan n\ran^{2k}}{\lan\tau+n^2\ran^{2a}} \sum\limits_{n_1}\int d\tau_1
\frac{\chi_{R_0}}{ \lan\tau_1+n_1^2\ran^{2b} \lan\tau_2\ran^{2c} \lan n_1\ran^{2k} \lan n_2\ran^{2s}}
\end{equation}
\begin{equation}\label{e.w2}
W_2:=\sup\limits_{n_1,\tau_1} \frac{1}{\lan n_1\ran^{2k} \lan\tau_1+n_1^2\ran^{2b}} \sum\limits_{n}\int d\tau
\frac{\lan n\ran^{2k}\chi_{R_1}}{ \lan\tau+n^2\ran^{2a} \lan\tau_2\ran^{2c} \lan n_2\ran^{2s}}
\end{equation}
\begin{equation}\label{e.w3}
W_3:=\sup\limits_{n_2,\tau_2} \frac{1}{\lan n_2\ran^{2s} \lan\tau_2\ran^{2c}} \sum\limits_{n}\int d\tau
\frac{\lan n\ran^{2k}\chi_{R_2}}{ \lan\tau+n^2\ran^{2a} \lan\tau_1+n_1^2\ran^{2b} \lan n_1\ran^{2k}}
\end{equation}
For later use, we recall that the dispersive relation of this bilinear estimate is:
\begin{equation}\label{e.dispersive-uv-periodic}
\tau+n^2-(\tau_1+n_1^2)-\tau_2 = n^2-n_1^2
\end{equation}
In order to define the regions $R_0, R_1, R_2$, we introduce the subsets:
\begin{equation}
\begin{split}
&A:=\left\{(n,n_1,\tau,\tau_1)\in\mathbb{Z}^2\times \mathbb{R}^2: |n|\lesssim 1\right\}, \\
&B:=\left\{(n,n_1,\tau,\tau_1)\in\mathbb{Z}^2\times \mathbb{R}^2: |n|\gg 1 \textrm{ and } |n|\sim |n_1|\right\}, \\
&C:=\left\{(n,n_1,\tau,\tau_1)\in\mathbb{Z}^2\times \mathbb{R}^2:
|n|\gg 1, |n|\nsim |n_1| \textrm{ and } |\tau+n^2|=L_{\max}\right\},
\end{split}
\end{equation}
where $L_{\max}:=\max\{|\tau+n^2|, |\tau_1+n_1^2|, |\tau_2|\}$. For
later reference, we denote also $N_{\max}:=\max\{|n|, |n_1|,
|n_2|\}$. Then, we put $R_0:=A\cup B\cup C$ and
\begin{equation}
\begin{split}
&R_1:=\left\{(n,n_1,\tau,\tau_1)\in\mathbb{Z}^2\times \mathbb{R}^2: |n|\gg 1, |n|\nsim |n_1| \textrm{ and }
|\tau_1+n_1^2|= L_{\max}\right\}, \\
&R_2:=\left\{(n,n_1,\tau,\tau_1)\in\mathbb{Z}^2\times \mathbb{R}^2: |n|\gg 1, |n|\nsim |n_1| \textrm{ and }
|\tau_2|= L_{\max}\right\}.
\end{split}
\end{equation}
We begin with the analysis of (\ref{e.w1}). In the region $A$, since $|n|\lesssim 1$, we have
\begin{equation*}
\begin{split}
&\sup\limits_{n,\tau} \frac{\lan n\ran^{2k}}{\lan\tau+n^2\ran^{2a}} \sum\limits_{n_1}\int d\tau_1
\frac{\chi_A}{ \lan\tau_1+n_1^2\ran^{2b} \lan\tau_2\ran^{2c} \lan n_1\ran^{2k} \lan n_2\ran^{2s}} \\
&\lesssim \sup\limits_{n,\tau} \frac{1}{\lan\tau+n^2\ran^{2a}}
\sum\limits_{n_1}\int d\tau_1 \frac{1}{ \lan\tau_1+n_1^2\ran^{2b}
\lan\tau_2\ran^{2c} \lan n_1\ran^{2k} \lan n_2\ran^{2s}}\\
&\lesssim \sup\limits_{\tau} \sum\limits_{n_1} \frac{1}{\lan\tau+n_1^2\ran^{2b+2c-1-}} \\
&\lesssim 1,
\end{split}
\end{equation*}
because $k,s\geq 0$, $a>0$ and $2b+2c>3/2$.

In the region $B$, we have $|n|\sim |n_1|$. Thus,
\begin{equation*}
\begin{split}
&\sup\limits_{n,\tau} \frac{\lan n\ran^{2k}}{\lan\tau+n^2\ran^{2a}}
\sum\limits_{n_1}\int d\tau_1
\frac{\chi_B}{ \lan\tau_1+n_1^2\ran^{2b} \lan\tau_2\ran^{2c} \lan n_1\ran^{2k} \lan n_2\ran^{2s}} \\
&\lesssim \sup\limits_{n,\tau} \frac{1}{\lan\tau+n^2\ran^{2a}}
\sum\limits_{n_1}\int d\tau_1 \frac{1}{ \lan\tau_1+n_1^2\ran^{2b}
\lan\tau_2\ran^{2c} \lan n_2\ran^{2s}}\\
&\lesssim \sup\limits_{\tau} \sum\limits_{n_1} \frac{1}{\lan\tau+n_1^2\ran^{2b+2c-1-}} \\
&\lesssim 1,
\end{split}
\end{equation*}
because $k,s\geq 0$, $a>0$ and $2b+2c>3/2$.

In the region $C$, we know that $|\tau+n^2|=L_{\max}$, $|n|\nsim
|n_1|$ and $|n|\gg 1$. Hence, the dispersive
relation~(\ref{e.dispersive-uv-periodic}) says that
$|\tau+n^2|=L_{\max}\gtrsim |n^2-n_1^2|\gtrsim N_{\max}^2$.
Therefore,
\begin{equation*}
\begin{split}
&\sup\limits_{n,\tau} \frac{\lan n\ran^{2k}}{\lan\tau+n^2\ran^{2a}}
\sum\limits_{n_1}\int d\tau_1
\frac{\chi_C}{ \lan\tau_1+n_1^2\ran^{2b} \lan\tau_2\ran^{2c} \lan n_1\ran^{2k} \lan n_2\ran^{2s}} \\
&\lesssim \sup\limits_{n,\tau} \frac{\lan
N_{\max}\ran^{2k-2s}}{\lan\tau+n^2\ran^{2a}} \sum\limits_{n_1}\int
d\tau_1 \frac{1}{ \lan\tau_1+n_1^2\ran^{2b}
\lan\tau_2\ran^{2c} }\\
&\lesssim \sup\limits_{n,\tau} \frac{\lan
N_{\max}\ran^{2k-2s}}{\lan N_{\max}\ran^{4a}}\sum\limits_{n_1} \frac{1}{\lan\tau+n_1^2\ran^{2b+2c-1-}} \\
&\lesssim 1,
\end{split}
\end{equation*}
since $k,s\geq 0$, $k<s+1$, $a=1/2-$ and $2b+2c>3/2$.

Putting together the estimates above, we conclude the desired
boundedness of~(\ref{e.w1}):
$$|W_1|\lesssim 1.$$

Next we estimate the contribution of~(\ref{e.w2}). In the region
$R_1$, we know that $|n|\gg 1$, $|n|\nsim |n_1|$ and
$|\tau_1+n_1^2|=L_{\max}$. So, the dispersive
relation~(\ref{e.dispersive-uv-periodic}) implies that
$|\tau_1+n_1^2|\gtrsim N_{\max}^2$. Thus,
\begin{equation*}
\begin{split}
&W_2:=\sup\limits_{n_1,\tau_1} \frac{1}{\lan n_1\ran^{2k}
\lan\tau_1+n_1^2\ran^{2b}} \sum\limits_{n}\int d\tau \frac{\lan
n\ran^{2k}\chi_{R_1}}{ \lan\tau+n^2\ran^{2a} \lan\tau_2\ran^{2c}
\lan n_2\ran^{2s}}\\
&\lesssim\sup\limits_{\tau_1} \sum\limits_{n}\int d\tau \frac{\lan
N_{\max}\ran^{2k-2s-4b}}{ \lan\tau+n^2\ran^{2a} \lan\tau_2\ran^{2c}}
\lesssim\sup\limits_{\tau_1} \sum\limits_{n}
\frac{1}{\lan\tau_1+n^2\ran^{2a+2c-1-}} \\
&\lesssim 1,
\end{split}
\end{equation*}
since $k,s\geq 0$, $k<s+1$ and $b=1/2-$, $2a+2c>3/2$.

Finally, we bound~(\ref{e.w3}) by noting that, in the region $R_2$,
it holds $|n|\gg 1$, $|n|\nsim |n_1|$ and $|\tau_2|=L_{\max}$. In
particular, the dispersive relation~(\ref{e.dispersive-uv-periodic})
forces $|\tau_2|\gtrsim N_{\max}^2$. This allows to obtain
\begin{equation*}
\begin{split}
&W_3:=\sup\limits_{n_2,\tau_2} \frac{1}{\lan n_2\ran^{2s}
\lan\tau_2\ran^{2c}} \sum\limits_{n}\int d\tau \frac{\lan
n\ran^{2k}\chi_{R_2}}{ \lan\tau+n^2\ran^{2a}
\lan\tau_1+n_1^2\ran^{2b} \lan n_1\ran^{2k}}\\
&\lesssim\sup\limits_{n_2,\tau_2} \sum\limits_{n}\int d\tau
\frac{\lan N_{\max}\ran^{2k-2s-4c}}{ \lan\tau+n^2\ran^{2a}
\lan\tau_1+n_1^2\ran^{2b} }\lesssim\sup\limits_{n_2,\tau_2}
\sum\limits_{n} \frac{1}{
\lan\tau_2+n_2(n_2-2n)\ran^{2a+2b-1-}}\\
&\lesssim 1,
\end{split}
\end{equation*}
since $k,s\geq 0$, $k<s+1$ and $2a+2b>3/2$.

This completes the proof of the proposition.
\end{proof}

\begin{proposition}\label{proposition-u2-periodic}
The bilinear estimate
\begin{equation}\label{Estimativa-u2}
\|u\overline{w}\|_{H^{-\frac{1}{2}+}_t
H_x^{s}}\lesssim
\|u\|_{X^{k,\frac{1}{2}+}}\|w\|_{X^{k,\frac{1}{2}+}}
\end{equation}
holds if $0\leq s\leq 2k$ and \small{$|k-s| < 1$}.
\end{proposition}
\begin{proof}
Similarly to the previous proposition, the relevant dispersive
relation is
\begin{equation}\label{e.dipersive-|u|2-periodic}
\tau - (\tau_1+n_1^2) - (\tau_2-n_2^2) = n_2^2-n_1^2
\end{equation}
and it suffices to bound the following contributions:
\begin{equation}\label{e.z1-periodic}
Z_1:=\sup\limits_{n,\tau} \frac{\lan n\ran^{2s}}{\lan\tau\ran^{2a}}
\sum\limits_{n_1}\int d\tau_1 \frac{\chi_{S_0}}{
\lan\tau_1+n_1^2\ran^{2b} \lan\tau_2-n_2^2\ran^{2c} \lan
n_1\ran^{2k} \lan n_2\ran^{2k}}
\end{equation}
\begin{equation}\label{e.z2}
Z_2:=\sup\limits_{n_1,\tau_1} \frac{1}{\lan n_1\ran^{2k}
\lan\tau_1+n_1^2\ran^{2b}} \sum\limits_{n}\int d\tau \frac{\lan
n\ran^{2s}\chi_{S_1}}{ \lan\tau\ran^{2a} \lan\tau_2-n_2^2\ran^{2c}
\lan n_2\ran^{2k}}
\end{equation}
\begin{equation}\label{e.z3}
Z_3:=\sup\limits_{n_2,\tau_2} \frac{1}{\lan n_2\ran^{2k}
\lan\tau_2-n_2^2\ran^{2c}} \sum\limits_{n}\int d\tau \frac{\lan
n\ran^{2s}\chi_{S_2}}{ \lan\tau\ran^{2a} \lan\tau_1+n_1^2\ran^{2b}
\lan n_1\ran^{2k}}
\end{equation}
where $S_0\cup S_1\cup S_2 = \mathbb{Z}^2\times\mathbb{R}^2$. To
define the regions $S_j$, $j=0,1,2$, we introduce the sets
\begin{equation}
\begin{split}
&E:=\left\{(n,n_1,\tau,\tau_1)\in\mathbb{Z}^2\times \mathbb{R}^2: |n|\lesssim 1\right\}, \\
&F:=\left\{(n,n_1,\tau,\tau_1)\in\mathbb{Z}^2\times \mathbb{R}^2: |n|\gg 1 \textrm{ and } |n_1|\sim |n_2|\right\}, \\
&G:=\left\{(n,n_1,\tau,\tau_1)\in\mathbb{Z}^2\times \mathbb{R}^2:
|n|\gg 1, |n_1|\nsim |n_2| \textrm{ and }
|\tau_1+n_1^2|=L_{\max}\right\},
\end{split}
\end{equation}
We put $S_1:=E\cup F\cup G$ and
\begin{equation}
\begin{split}
&S_0:=\left\{(n,n_1,\tau,\tau_1)\in\mathbb{Z}^2\times \mathbb{R}^2:
|n|\gg 1, |n_1|\nsim |n_2| \textrm{ and }
|\tau|= L_{\max}\right\}, \\
&S_2:=\left\{(n,n_1,\tau,\tau_1)\in\mathbb{Z}^2\times \mathbb{R}^2:
|n|\gg 1, |n_1|\nsim |n_2| \textrm{ and } |\tau_2-n_2^2|=
L_{\max}\right\}.
\end{split}
\end{equation}
We can estimate $Z_1$ as follows. In the region $S_0$, since
$|n_1|\nsim |n_2|$, we have either $|n_1|\gg |n_2|$ or $|n_2|\gg
|n_1|$. By symmetry reasons, we can suppose that, without loss of
generality, $|n_2|\gg |n_1|$. In this case, $|\tau|\gtrsim n_2^2$
and $|n|\sim |n_2|$. So,
\begin{equation}
\begin{split}
&Z_1:=\sup\limits_{n,\tau} \frac{\lan n\ran^{2s}}{\lan\tau\ran^{2a}}
\sum\limits_{n_1}\int d\tau_1 \frac{\chi_{S_0}}{
\lan\tau_1+n_1^2\ran^{2b} \lan\tau_2-n_2^2\ran^{2c} \lan
n_1\ran^{2k} \lan n_2\ran^{2k}}\\
&\lesssim\sup\limits_{n,\tau} \frac{\lan
n\ran^{2s-2k}}{\lan\tau\ran^{2a}} \sum\limits_{n_1}\int d\tau_1
\frac{\chi_{S_0}}{ \lan\tau_1+n_1^2\ran^{2b}
\lan\tau_2-n_2^2\ran^{2c}}\\
&\lesssim\sup\limits_{n,\tau} \sum\limits_{n_1} \frac{1}{
\lan\tau+n_1^2-n_2^2\ran^{2b+2c-1-} \lan\tau_2-n_2^2\ran^{2c} \lan
n_1\ran^{2k} \lan n_2\ran^{2k}}\\
&\lesssim 1,
\end{split}
\end{equation}
since $k,s\geq 0$, $s<k+1$ and $2b+2c>3/2$.

Now we will bound the expression $Z_2$. In the region $E$, it holds
$|n|\lesssim 1$. Hence,
\begin{equation}
\begin{split}
&\sup\limits_{n_1,\tau_1} \frac{1}{\lan n_1\ran^{2k}
\lan\tau_1+n_1^2\ran^{2b}} \sum\limits_{n}\int d\tau \frac{\lan
n\ran^{2s}\chi_{E}}{ \lan\tau\ran^{2a} \lan\tau_2-n_2^2\ran^{2c}
\lan n_2\ran^{2k}}\\
&\lesssim\sup\limits_{n_1,\tau_1} \sum\limits_{|n|\lesssim 1}\int
d\tau \frac{1}{ \lan\tau\ran^{2a} \lan\tau_2-n_2^2\ran^{2c}}
\lesssim\sup\limits_{n_1,\tau_1} \sum\limits_{|n|\lesssim 1}
\frac{1}{ \lan\tau_1+n_2^2\ran^{2a+2c-1-}}\\
&\lesssim 1.
\end{split}
\end{equation}
In the region $F$, we get $|n_1|\sim |n_2|$ so that
\begin{equation}
\begin{split}
&\sup\limits_{n_1,\tau_1} \frac{1}{\lan n_1\ran^{2k}
\lan\tau_1+n_1^2\ran^{2b}} \sum\limits_{n}\int d\tau \frac{\lan
n\ran^{2s}\chi_{F}}{ \lan\tau\ran^{2a} \lan\tau_2-n_2^2\ran^{2c}
\lan n_2\ran^{2k}}\\
&\lesssim\sup\limits_{n_1,\tau_1} \frac{1}{
\lan\tau_1+n_1^2\ran^{2b}} \sum\limits_{n}\int d\tau \frac{\lan
n\ran^{2s-4k}}{ \lan\tau\ran^{2a} \lan\tau_2-n_2^2\ran^{2c}}\\
&\lesssim\sup\limits_{n_1,\tau_1} \sum\limits_{n}
\frac{1}{ \lan\tau_1 + (n-n_1)^2\ran^{2a+2c-1-}}\\
&\lesssim 1,
\end{split}
\end{equation}
because $0\leq s\leq 2k$ and $2a+2c>3/2$. In the region $G$, the
dispersive relation~(\ref{e.dipersive-|u|2-periodic}) combined with
the assumptions $|n|\gg 1$, $|n_1|\nsim |n_2|$ and $|\tau_1+n_1^2|=
L_{\max}$ implies that $|\tau_1+n_1^2|\gtrsim N_{\max}^2$. Without
loss of generality, we can suppose that $|n_1|\ll |n_2|$. Then,
\begin{equation}
\begin{split}
&\sup\limits_{n_1,\tau_1} \frac{1}{\lan n_1\ran^{2k}
\lan\tau_1+n_1^2\ran^{2b}} \sum\limits_{n}\int d\tau \frac{\lan
n\ran^{2s}\chi_{G}}{ \lan\tau\ran^{2a} \lan\tau_2-n_2^2\ran^{2c}
\lan n_2\ran^{2k}}\\
&\lesssim\sup\limits_{n_1,\tau_1} \sum\limits_{n}\int d\tau
\frac{\lan n\ran^{2s-2k-4b}}{ \lan\tau\ran^{2a}
\lan\tau_2-n_2^2\ran^{2c}}\\
&\lesssim 1,
\end{split}
\end{equation}
since $0\leq k,s$ and $s<k+1$, $2a+2c>3/2$. Collecting these
estimates, we conclude
\begin{equation}
|Z_2|\lesssim 1.
\end{equation}
Finally, the expression~(\ref{e.z3}) can be controlled if we notice
that $|n|\gg 1$, $|n_1|\nsim |n_2|$ and $|\tau_2-n_2^2|= L_{\max}$
implies $|\tau_2-n_2^2|\gtrsim N_{\max}^2$. In particular,
\begin{equation}
\begin{split}
&Z_3:=\sup\limits_{n_2,\tau_2} \frac{1}{\lan n_2\ran^{2k}
\lan\tau_2-n_2^2\ran^{2c}} \sum\limits_{n}\int d\tau \frac{\lan
n\ran^{2s}\chi_{S_2}}{ \lan\tau\ran^{2a} \lan\tau_1+n_1^2\ran^{2b}
\lan n_1\ran^{2k}}\\
&\lesssim\sup\limits_{n_2,\tau_2} \sum\limits_{n}\int d\tau
\frac{\lan N_{\max}\ran^{2s-2k-4c}}{ \lan\tau\ran^{2a}
\lan\tau_1+n_1^2\ran^{2b}}\\
&\lesssim\sup\limits_{n_2,\tau_2} \sum\limits_{n} \frac{1}{
\lan\tau_2-(n-n_2)^2\ran^{2a+2b-1-}}\\
&\lesssim 1,
\end{split}
\end{equation}
whenever $k,s\geq 0$, $s<k+1$ and $2a+2b>3/2$.

This finishes the proof of the proposition.
\end{proof}

\begin{remark}Again, once the bilinear estimates in propositions~\ref{proposition-uv-periodic} and \ref{proposition-u2-periodic} are proved, one can show the theorem \ref{local-theorem-periodic} by standard arguments (e.g., see the works \cite{KPV}, \cite{Bekiranov1} and \cite{Ginibre}).
\end{remark}

\begin{remark}After the completion of this work, Angulo, Corcho and Hakkaev~\cite{ACH} improved the bilinear estimate of proposition~\ref{proposition-uv-periodic} (for the coupling term $uv$) so that we can include the boundary case $|k-s|=1$ in the statement of our proposition (if one is willing to modify a little bit the definition of the Bourgain spaces). Nevertheless, it is possible to show that the same method leads to an improved bilinear estimate of proposition~\ref{proposition-u2-periodic} (for the coupling term $|u|^2$) in order to include again the boundary case $|k-s|=1$. Hence, it follows that the local well-posedness result of theorem~\ref{local-theorem-periodic} holds for any pair of indices $(k,s)$ verifying $0\leq s\leq 2k$ and $|k-s|\leq 1$.
\end{remark}


\subsection{Counter-Examples II: the periodic case}  The next results prove that the
bilinear estimates derived in propositions
\ref{proposition-uv-periodic} and \ref{proposition-u2-periodic} are sharp.

\begin{proposition}For any $b_1,b_2\in\mathbb{R}$, the estimate $\|uv\|_{X^{k,-\frac{1}{2}+}}\lesssim
\|u\|_{X^{k,b_1}} \|v\|_{H_t^{b_2}H_x^{s}}$ holds only when $s\geq 0$ and $k < s+1$.
\end{proposition}

\begin{proof} Firstly, we fix $N\gg 1$ a large integer and define
\begin{displaymath}
a_n= \left\{ \begin{array}{ll}
1 & \textrm{if $n=N$}\\
0 & \textrm{otherwise}
\end{array} \right.
\end{displaymath}
and
\begin{displaymath}
b_n= \left\{ \begin{array}{ll}
1 & \textrm{if $n=-2N$}\\
0 & \textrm{otherwise}
\end{array} \right.
\end{displaymath}

Let $f$ and $g$ be given by
$\hat{f}(n,\tau)=a_n\chi_{[-1,1]}(\tau+n^2)$\; and\;
$\hat{g}(n,\tau)=b_n\chi_{[-1,1]}(\tau)$. Taking
into account the dispersive relation
$\tau+n^2-(\tau_1+n_1^2)-\tau_2 = n^2-n_1^2$, we can easily
compute that
\begin{equation*}
\|fg\|_{X^{k,-1/2+}} \simeq N^k,\quad \|f\|_{X^{k,b_1}} \simeq
N^k\quad  \text{and}\quad \|g\|_{H_t^{b_2}H_x^s} \simeq N^s
\end{equation*}
Hence, the bound $\|fg\|_{X^{k,-\frac{1}{2}+}}\lesssim
\|f\|_{X^{k,\frac{1}{2}}} \|g\|_{H_t^{\frac{1}{2}} H_x^{s}}$
implies $N^k\lesssim N^{k+s}$, consequently, $s\geq 0$.

Secondly, define
\begin{displaymath}
d_n= \left\{ \begin{array}{ll}
1 & \textrm{if $n=N$}\\
0 & \textrm{otherwise}
\end{array} \right.
\end{displaymath}
and
\begin{displaymath}
c_n= \left\{ \begin{array}{ll}
1 & \textrm{if $n=0$}\\
0 & \textrm{otherwise}
\end{array} \right.
\end{displaymath}

Let $p$ and $q$ be $\hat{p}(n,\tau)=c_n\chi_1(\tau+n^2)$ and
$\hat{q}(n,\tau)=d_n\chi_1(\tau)$. Again, it is not hard to see
that
\begin{equation*}
\|pq\|_{X^{k,-1/2+}} \simeq \frac{N^k}{N^{1-}}
\end{equation*}
\begin{equation*}
\|p\|_{X^{k,b_1}} \simeq 1
\end{equation*}
\begin{equation*}
\|q\|_{H_t^{b_2}H_x^s} \simeq N^s
\end{equation*}
Hence, the bound $\|pq\|_{X^{k,-\frac{1}{2}+}}\lesssim
\|p\|_{X^{k,\frac{1}{2}}} \|q\|_{H_t^{\frac{1}{2}} H_x^{s}}$
implies $\frac{N^k}{N^{1-}}\lesssim N^{s}$, i.e., $k < s+1$.
\end{proof}

\begin{proposition}For any $b_1,b_2\in\mathbb{R}$, the estimate $\|u_1 \overline{u_2}\|_{H_t^{-\frac{1}{2}+} H_x^s}\lesssim \|u_1\|_{X^{k,b_1}} \|u_2\|_{X^{k,b_2}}$ holds only if
$s\leq 2k$ and $s < k+1$.
\end{proposition}

\begin{proof}For a fixed large integer $N\gg 1$, define
\begin{displaymath}
a_n= \left\{ \begin{array}{ll}
1 & \textrm{if $n=N$}\\
0 & \textrm{otherwise}
\end{array} \right.
\end{displaymath}
\begin{displaymath}
b_n= \left\{ \begin{array}{ll}
1 & \textrm{if $n=-N-1$}\\
0 & \textrm{otherwise}
\end{array} \right.
\end{displaymath}
\begin{displaymath}
c_n= \left\{ \begin{array}{ll}
1 & \textrm{if $n=0$}\\
0 & \textrm{otherwise}
\end{array} \right.
\end{displaymath}
\begin{displaymath}
d_n= \left\{ \begin{array}{ll}
1 & \textrm{if $n=N$}\\
0 & \textrm{otherwise}
\end{array} \right.
\end{displaymath}
Putting $\hat{f_1}(n,\tau) = a_n \chi_1(\tau+n^2)$,
$\hat{f_2}(n,\tau) = b_n \chi_1(\tau+n^2)$ and $\hat{g_1}(n,\tau)
= c_n \chi_1(\tau+n^2)$, $\hat{g_2}(n,\tau) = d_n
\chi_1(\tau+n^2)$, a simple calculation (based on the dispersive
relation $\tau - (\tau_1+n_1^2) - (\tau_2+n_2^2) = n_2^2 - n_1^2$)
gives that
\begin{equation*}
\|f_1\overline{f_2}\|_{H_t^{-1/2+} H_x^s}\simeq N^s \quad
\textrm{and} \quad \|g_1 \overline{g_2}\|_{H_t^{-1/2+}
H_x^s}\simeq N^{s-1+},
\end{equation*}
\begin{equation*}
\|f_1\|_{X^{k,b_1}}\simeq N^k \quad \textrm{and} \quad
\|g_1\|_{X^{k,b_1}}\simeq 1,
\end{equation*}
\begin{equation*}
\|f_2\|_{X^{k,b_2}}\simeq N^k\simeq \|g_2\|_{X^{k,b_2}}.
\end{equation*}
Therefore, the bound $\|u_1 \overline{u_2}\|_{H_t^{-\frac{1}{2}+}
H_x^s}\lesssim \|u_1\|_{X^{k,b_1}}\|u_2\|_{X^{k,b_2}}$ says that $N^s\lesssim N^{2k}$ and
$N^{s-1+}\lesssim N^k$, i.e., $s\leq 2k$ and $s< k+1$.
\end{proof}



\section{Global well-posedness below $L^2\times
L^2$}\label{section-global-results}

This section is devoted to the proof of the global well-posedness result stated in theorem~\ref{global-theorem} via the I-method of Colliander, Keel, Staffilani, Takaoka and Tao.

\subsection{The I-operator} Let $m(\xi)$ be a smooth non-negative symbol on $\mathbb{R}$ which equals $1$ for $|\xi|\leq 1$ and equals $|\xi|^{-1}$ for $|\xi|\geq 2$. For any $N\geq 1$ and $\alpha\in\mathbb{R}$, denote by $I_N^{\alpha}$ the Fourier multiplier
\begin{equation*}
\widehat{I_N^{\alpha} f}(\xi) = m\left(\frac{\xi}{N}\right)^{\alpha} \widehat{f}(\xi).
\end{equation*}

We recall the following abstract interpolation lemma:

\begin{lemma}[Lemma 12.1 of~\cite{CKSTT2}]\label{l.interpolation} Let $\alpha_0>0$ and $n\geq 1$. Suppose $Z, X_1, \dots, X_n$ are translation-invariant Banach spaces and $T$ is a translation invariant $n$-linear operator such that
\begin{equation*}
\|I_1^{\alpha} T(u_1,\dots,u_n)\|_Z\lesssim \prod\limits_{j=1}^n \|I_1^{\alpha}u_j\|_{X_j},
\end{equation*}
for all $u_1,\dots,u_n$ and $0\leq\alpha\leq\alpha_0$. Then,
\begin{equation*}
\|I_N^{\alpha}T(u_1,\dots,u_n)\|_Z\lesssim \prod\limits_{j=1}^n \|I_N^{\alpha}u_j\|_{X_j},
\end{equation*}
for all $u_1,\dots,u_n$, $0\leq\alpha\leq\alpha_0$ and $N\geq 1$. Here, the implied constant is independent of $N$.
\end{lemma}

After these preliminaries, we are ready to show a variant of the local well-posedness theorem~\ref{local-theorem-continuous}.

\subsection{Local well-posedness revisited}

In the sequel, we take $N\gg 1$ a large integer and we denote by $I$ the operator $I:= I_N^{-s}$ for a given $s\in\mathbb{R}$.

\begin{proposition}\label{p.local-I}For all $(u_0,v_0)\in H^s(\mathbb{R})\times H^s(\mathbb{R})$ and $s\geq -1/4$, the Schr\"odinger-Debye system (\ref{S-Debye}) has a unique local-in-time solution $(u(t), v(t))$ defined on the time interval $[0,\delta]$ for some $\delta\leq 1$ satisfying
\begin{equation}\label{e.local-I}
\delta\sim (\|Iu_0\|_{L_x^2}+\|Iv_0\|_{L_x^2})^{-4/3-}.
\end{equation}
Furthermore, $\|Iu\|_{X^{0,1/2+}}\lesssim\|Iu_0\|_{L^2}$ and
$\|Iv\|_{X^{0,1/2+}} \lesssim \|Iu_0\|_{L^2}+\|Iv_0\|_{L^2}$.
\end{proposition}

\begin{proof}Applying the I-operator to the Schr\"odinger-Debye system (\ref{S-Debye}), we get
\begin{equation}
\label{I-SD}
\begin{cases}
i\p_tIu+ \tfrac{1}{2}\partial_x^2Iu=I(uv), \\
\sigma \p_tIv + Iv = \epsilon I(|u|^2),\\
u(x,0)=u_0(x) , \quad v(x,0)=v_0(x).
\end{cases}
\end{equation}
To solve this problem, we denote by $\Phi_1(Iu,Iv)$ and
$\Phi_2(Iu,Iv)$ the integral maps associated to this system, so that
our task is to find a fixed point of $(\Phi_1,\Phi_2)$. To
accomplish this objective, note that, by standard arguments, the lemma~\ref{l.time}, the
interpolation lemma~\ref{l.interpolation} combined with the bilinear
estimates in the corollaries~\ref{corollary-uv-continuous}
and~\ref{corollary-|u2|-continuous} give the estimates

\begin{equation*}
\begin{split}
&\|\Phi_1(Iu,Iv)\|_{X^{0,1/2+}}\leq C\|I u_0\|_{L^2_x} + C\delta^{3/4-}\|I
u\|_{X^{0,1/2+}}\|I v\|_{H_t^{1/2+}L_x^2},\\
&\|\Phi_2(Iu, Iv)\|_{H_t^{1/2+}L_x^2}\leq C\|I v_0\|_{L_x^2} +
C\delta^{3/4-}\|Iu\|_{X^{0,1/2+}}^2,
\end{split}
\end{equation*}
where $Iu,Iv\in X^{0,1/2+}$ are defined in the interval $[0,\delta]$.

Taking $R_1=2C\|Iu_0\|_{L_x^2}$ and $R_2=2C(\|Iu_0\|_{L_x^2} +
\|Iv_0\|_{L_x^2})$, we conclude that $(\Phi_1,\Phi_2)$ has an unique
fixed point $(Iu,Iv)$ on the product $B(R_1)\times B(R_2)$ of balls
of radii $R_1$ and $R_2$. Moreover,
\begin{equation*}
\delta\sim (\|Iu_0\|_{L_x^2}+\|Iv_0\|_{L_x^2})^{-4/3-}
\end{equation*}
This completes the proof of the proposition.
\end{proof}

Once a local well-posedness result for the modified system $(\ref{I-SD})$ was obtained, we will study the behavior of the $L^2$-conservation law under the $I$-operator.

\subsection{Modified energy}

We consider the modified energy $E(Iu) = \|Iu\|_{L_x^2}^2$. Note
that, since $(Iu,Iv)$ verify the system (\ref{I-SD}), we have
\begin{equation*}
\begin{split}
\frac{d}{dt} E(Iu)(t) &= \int \partial_t Iu \cdot I\overline{u} +
\int Iu \cdot \partial_t I\overline{u} \\
&= -\frac{1}{i}\int\partial_x^2 Iu\cdot I\overline{u} +
\frac{1}{i}\int I(uv) I\overline{u} +\frac{1}{i}\int Iu
\partial_x^2I\overline{u} -\frac{1}{i}\int Iu
I(\overline{uv}) \\
&= \frac{1}{i}\int \partial_x Iu\cdot\partial_x I\overline{u} +
\frac{1}{i}\int \left(I(uv)-Iu Iv\right)I\overline{u} +
\frac{1}{i}\int Iu Iv I\overline{u} \\
&-\frac{1}{i}\int \partial_x Iu \cdot\partial_x I\overline{u}
-\frac{1}{i}\int Iu \overline{(I(uv)-Iu Iv)}-\frac{1}{i}\int Iu
I\overline{u} Iv \\
&= 2\Im \int \left(I(uv)-IuIv\right)I\overline{u}.
\end{split}
\end{equation*}

Now we are going to see that this formula leads naturally to an \emph{almost conservation law}.

\subsection{Almost conservation of the modified energy}

For later use, we need the following refined Strichartz estimate:

\begin{lemma}\label{l.bourgain}We have
\begin{equation*}
\|(D_x^{1/2}f)\cdot g\|_{L_{xt}^2}\lesssim \|f\|_{X^{0,1/2+}} \|g\|_{X^{0,1/2+}},
\end{equation*}
if $|\xi_1|\gg |\xi_2|$ for any $|\xi_1|\in \textrm{supp}(\widehat{f}), |\xi_2|\in \textrm{supp}(\widehat{g})$. Moreover, this estimate is true if $f$ and/or $g$ is replaced by its complex conjugate in the left-hand side of the inequality.
\end{lemma}

\begin{proof}See lemma 7.1 of~\cite{CKSTT3} or lemma 4.2 of~\cite{Grunrock}.
\end{proof}

\begin{lemma}\label{l.ac}For $s>-1/4$ and any parameter $1/8<\ell<1/4$, it holds
\begin{equation*}
|E(Iu)(\delta) - E(Iu)(0)|\lesssim
N^{-2\ell+}\delta^{\tfrac{1}{2}-2\ell-}\|Iu\|_{X^{0,1/2+}}^2
\|Iv\|_{H_t^{1/2+}L_x^2}
\end{equation*}
\end{lemma}

\begin{proof} Since we already know that
\begin{equation*}
E(Iu)(\delta) - E(Iu)(0) = \int_0^{\delta}\frac{d}{dt}E(Iu)(t) dt
= 2\Im \int_0^{\delta}\int\left(I(uv)-IuIv\right)I\overline{u},
\end{equation*}
it suffices to show that
\begin{equation}\label{e.ac}
\left|\int_0^{\delta}\int\left(I(uv)-IuIv\right)I\overline{u}\right|\lesssim
N^{-1/2+}\delta^{1/4-}\|Iu\|_{X^{0,1/2}}^2\|Iv\|_{H_t^{1/2+}L_x^2}.
\end{equation}
By Parseval, our task is to prove that
\begin{equation*}
\begin{split}
I:=\int_0^{\delta}\int_{\xi_1+\xi_2-\xi_3=0}
&\left|\frac{m(\xi_1+\xi_2)-m(\xi_1)m(\xi_2)}{m(\xi_1)m(\xi_2)}\right|\widehat{u}(\xi_1,t)\widehat{v}(\xi_2,t)\widehat{\overline{w}}(\xi_3,t)
\\ &\lesssim
N^{-1/2+}\delta^{1/4-}\|u\|_{X^{0,1/2}}\|v\|_{H_t^{1/2+}L_x^2}\|w\|_{X^{0,1/2}}
\end{split}
\end{equation*}
We decompose the frequencies $\xi_j$, $j=1,2,3$ into dyadic blocks
$|\xi_j|\sim N_j$. Before starting the proof of this inequality,
we note that the multiplier
$M:=\frac{m(\xi_1+\xi_2)-m(\xi_1)m(\xi_2)}{m(\xi_1)m(\xi_2)}$
satisfies
\begin{itemize}
\item if $|\xi_1|\ll|\xi_2|, |\xi_1|\ll N$, then $$|M|\lesssim \left|\frac{m(\xi_1+\xi_2)-m(\xi_2)}{m(\xi_2)}\right|\lesssim
\left|\frac{\nabla m(\xi_2) \xi_1}{m(\xi_2)}\right|\lesssim
\frac{N_1}{N_2}$$.
\item similarly, if $|\xi_2|\ll|\xi_1|, |\xi_2|\ll N$, then $M\lesssim
N_2/N_1$.
\item if $|\xi_1|\ll|\xi_2|, |\xi_1|\gtrsim N$, then $$|M|\lesssim \frac{1}{m(\xi_1)}\lesssim \left(\frac{N_1}{N}\right)^{1/4-},$$
because $s>-1/4$.
\item similarly, if $|\xi_2|\ll|\xi_1|, |\xi_2|\gtrsim N$, then $|M|\lesssim
(N_2/N)^{1/4-}$.
\item finally, if $|\xi_1|\sim|\xi_2|\gtrsim N$, then
$$|M|\lesssim \frac{1}{m(\xi_1)m(\xi_2)}\lesssim
\left(\frac{N_1}{N}\right)^{1/2-}.$$
\end{itemize}
Therefore, we can bound $I$ as follows:
\begin{itemize}
\item When $|\xi_1|\ll|\xi_2|, |\xi_1|\ll N$, we have
$|\xi_3|\sim|\xi_2|\gg|\xi_1|$. Thus, from the
lemma~\ref{l.bourgain},
\begin{equation*}
\begin{split}
I&\lesssim \frac{N_1}{N_2} \frac{1}{N_3^{1/2}} \|D_x^{1/2}w \cdot
u\|_{L_{xt}^2}\|v\|_{L_{xt}^2} \\&\lesssim N^{-1/2+}\delta^{1/2}
N_{\max}^{0-}
\|u\|_{X^{0,1/2+}}\|v\|_{H_t^{1/2+}L_x^2}\|w\|_{X^{0,1/2+}}
\end{split}
\end{equation*}
\item if $|\xi_2|\ll|\xi_1|, |\xi_2|\ll N$, we also have
$|\xi_1|\sim|\xi_3|$; in this case, by duality, the lemma~\ref{l.time} and the bilinear
estimate of proposition~\ref{proposition-|u|2-continuous},
\begin{equation*}
\begin{split}
I&\lesssim \frac{N_2}{N_1}
\|u\overline{w}\|_{H_t^{-1/2-}L_x^2}\|v\|_{H_t^{1/2+}L_x^2} \\
&\lesssim \frac{N_2}{N_1}\delta^{\tfrac{1}{2}-2\ell-}
\|u\overline{w}\|_{H_t^{-2\ell-}L_x^2} \|v\|_{H_t^{1/2+}L_x^2}\\
&\lesssim
\delta^{\tfrac{1}{2}-2\ell-}\frac{N_2}{N_1}\|u\|_{X^{-\ell,1/2+}}\|w\|_{X^{-\ell,1/2+}}
\|v\|_{H_t^{1/2+}L_x^2} \\
&\lesssim
\delta^{\tfrac{1}{2}-2\ell-}\frac{N_2}{N_1}\frac{1}{N_1^{2\ell}}\|u\|_{X^{0,1/2+}}\|w\|_{X^{0,1/2+}}
\|v\|_{H_t^{1/2+}L_x^2} \\
&\lesssim N^{-2\ell+}\delta^{\tfrac{1}{2}-2\ell-} N_{\max}^{0-}\|u\|_{X^{0,1/2+}}
\|v\|_{H_t^{1/2+}L_x^2}\|w\|_{X^{0,1/2+}}.
\end{split}
\end{equation*}
\item when $|\xi_1|\ll|\xi_2|, N\lesssim |\xi_1|$, we know that
$|\xi_3|\sim|\xi_2|\gg|\xi_1|$, so that
\begin{equation*}
\begin{split}
I&\lesssim
\left(\frac{N_1}{N_2}\right)^{1/4-}\frac{1}{N_3^{1/2}}\|D_x^{1/2}w\cdot
u\|_{L^2}\|v\|_{L^2} \\
&\lesssim N^{-1/2+}\delta^{1/2}N_{\max}^{0-}\|u\|_{X^{0,1/2+}}
\|v\|_{H_t^{1/2+}L_x^2}\|w\|_{X^{0,1/2+}}.
\end{split}
\end{equation*}
\item if $|\xi_2|\ll |\xi_1|, N\lesssim |\xi_2|$, we have
$|\xi_1|\sim|\xi_3|$; thus,
\begin{equation*}
\begin{split}
I&\lesssim \left(\frac{N_2}{N}\right)^{1/4-}\|u\overline{w}\|_{H_t^{-1/2-}L_x^2}\|v\|_{H_t^{1/2+}L_x^2} \\
&\lesssim\left(\frac{N_2}{N}\right)^{1/4-}\delta^{\tfrac{1}{2}-2\ell-}\|u\overline{w}\|_{H_t^{-2\ell-}L_x^2} \|v\|_{H_t^{1/2+}L_x^2}\\
&\lesssim
\delta^{\tfrac{1}{2}-2\ell-}\left(\frac{N_2}{N}\right)^{1/4-}\|u\|_{X^{-\ell,1/2+}}\|w\|_{X^{-\ell,1/2+}}
\|v\|_{H_t^{1/2+}L_x^2} \\
&\lesssim
\delta^{\tfrac{1}{2}-2\ell-}\left(\frac{N_2}{N}\right)^{1/4-}\frac{1}{N_1^{2\ell}}\|u\|_{X^{0,1/2+}}\|w\|_{X^{0,1/2+}}
\|v\|_{H_t^{1/2+}L_x^2} \\
&\lesssim N^{-2\ell+}\delta^{\tfrac{1}{2}-2\ell-} N_{\max}^{0-}\|u\|_{X^{0,1/2+}}
\|v\|_{H_t^{1/2+}L_x^2}\|w\|_{X^{0,1/2+}}.
\end{split}
\end{equation*}
\item finally, when $|\xi_1|\sim|\xi_2|\gtrsim N$, we have two
possibilities: either $|\xi_1|\ll |\xi_3|$, so that
\begin{equation*}
\begin{split}
I&\lesssim \left(\frac{N_1}{N_2}\right)^{1/2-} \frac{1}{N_1^{1/2}}
\|D_x^{1/2}u \cdot w\|_{L_{xt}^2}\|v\|_{L_{xt}^2} \\&\lesssim
N^{-1/2+}\delta^{1/2} N_{\max}^{0-}
\|u\|_{X^{0,1/2+}}\|v\|_{H_t^{1/2+}L_x^2}\|w\|_{X^{0,1/2+}}
\end{split}
\end{equation*}
or $|\xi_1|\sim|\xi_3|$ implying
\begin{equation*}
\begin{split}
I&\lesssim \frac{N_2}{N_1}
\|u\overline{w}\|_{H_t^{-1/2-}L_x^2}\|v\|_{H_t^{1/2+}L_x^2} \\
&\lesssim \frac{N_2}{N_1}\delta^{\tfrac{1}{2}-2\ell-}
\|u\overline{w}\|_{H_t^{-2\ell-}L_x^2} \|v\|_{H_t^{1/2+}L_x^2}\\
&\lesssim
\delta^{\tfrac{1}{2}-2\ell-}\frac{N_2}{N_1}\|u\|_{X^{-\ell,1/2+}}\|w\|_{X^{-\ell,1/2+}}
\|v\|_{H_t^{1/2+}L_x^2} \\
&\lesssim
\delta^{\tfrac{1}{2}-2\ell-}\frac{N_2}{N_1}\frac{1}{N_1^{2\ell}}\|u\|_{X^{0,1/2+}}\|w\|_{X^{0,1/2+}}
\|v\|_{H_t^{1/2+}L_x^2} \\
&\lesssim N^{-2\ell+}\delta^{\tfrac{1}{2}-2\ell-} N_{\max}^{0-}\|u\|_{X^{0,1/2+}}
\|v\|_{H_t^{1/2+}L_x^2}\|w\|_{X^{0,1/2+}}.
\end{split}
\end{equation*}
\end{itemize}
Hence, in any case, we proved that
\begin{equation*}
I\lesssim N^{-2\ell+}\delta^{\tfrac{1}{2}-2\ell-} N_{\max}^{0-}\|u\|_{X^{0,1/2+}}
\|v\|_{H_t^{1/2+}L_x^2}\|w\|_{X^{0,1/2+}}.
\end{equation*}
Summing up over the dyadic blocks, we complete the proof of the
lemma.
\end{proof}

\subsection{Global existence}

Recall that $\|Iu_0\|_{L_x^2}\lesssim N^{-s} \|u_0\|_{H^s}$,
$\|Iv_0\|_{L_x^2}\lesssim N^{-s} \|v_0\|_{H^s}$,
$\|Iu\|_{X^{0,b}}\lesssim N^{-s} \|u\|_{X^{s,b}}$ and $\|Iv\|_{H_t^{1/2+}L^2_x}\lesssim N^{-s} \|v\|_{H_t^{1/2+}H^s_x}$. Applying the
local result of proposition~\ref{p.local-I}, we get the existence of
solutions on a time interval $[0,\delta]$, where $\delta\sim
N^{4s/3-}$. Also, they verify
\begin{equation*}
\|Iu\|_{X^{0,1/2+}}+\|Iv\|_{H_t^{1/2+}L_x^2}\lesssim N^{-s}.
\end{equation*}
By the lemma~\ref{l.ac}, for a given parameter $1/8<\ell<1/4$, we obtain
\begin{equation*}
|E(Iu)(\delta)-E(Iu)(0)|\lesssim N^{-2\ell+}\delta^{\tfrac{1}{2}-2\ell-}N^{-3s}.
\end{equation*}

On the other hand, using the lemma~\ref{l.time}, the bilinear estimate of
corollary~\ref{corollary-|u2|-continuous}, the interpolation
lemma~\ref{l.interpolation} and the local result of
proposition~\ref{p.local-I}, we get
\begin{equation*}
\begin{split}
\|Iv(\delta_0)\|_{L^2_x}&\leq e^{-\delta_0/\sigma}\|Iv_0\|_{L^2_x} +
\frac{1}{\sigma}\int_0^{\delta_0}e^{-(\delta_0-t)/\sigma}\|
I(|u(t)|^2)\|_{L^2_x} dt\\
&\leq e^{-\delta_0/\sigma}\|Iv_0\|_{L^2_x} +
C\left\|e^{-(\delta_0-t)/\sigma}\right\|_{H_t^{\tfrac{1}{2}-}}
\|I(|u|^2)\|_{H_t^{-\tfrac{1}{2}+}L_x^2}\\
&\leq e^{-\delta_0/\sigma}\|Iv_0\|_{L^2_x} +
C\left(\int_0^{\delta_0}e^{-2(\delta_0-t)/\sigma}dt\right)^{1/2} \delta_0^{\tfrac{1}{4}-}\|I(|u|^2)\|_{H_t^{-\tfrac{1}{4}-}L_x^2}\\
&\leq e^{-\delta_0/\sigma}\|Iv_0\|_{L^2_x} +
C\left(1-e^{-2\delta_0/\sigma}\right)^{1/2}\delta_0^{\tfrac{3}{4}-}\|Iu\|_{X^{0,\tfrac{1}{2}+}}^2\\
&\leq e^{-\delta_0/\sigma}\|Iv_0\|_{L^2_x} +
\left(1-e^{-\delta_0/\sigma}\right)\cdot C\delta_0^{\tfrac{1}{4}-}\|Iu_0\|_{L^2_x}^{2}\\
&\leq \max\{\|Iv_0\|_{L_x^2},C\delta_0^{\tfrac{1}{4}-}\|Iu_0\|_{L_x^2}^2\}\\
&\lesssim N^{-s}
\end{split}
\end{equation*}
for any $\delta_0\sim N^{4s-}$. In particular, since $\|Iu\|_{L_t^{\infty}L_x^2}\lesssim \|Iu\|_{X^{0,1/2+}}\lesssim N^{-s}$, we can iterate the previous estimate ($\delta/\delta_0$ times) to obtain $\|Iv(\delta)\|_{L_x^2}\lesssim N^{-s}$.

Finally, we observe that one can iterate the local result to cover the time interval
$[0,T]$ if these estimates hold after $T/\delta$ steps. In other
words, the existence of a solution on the time interval $[0,T]$ is
guaranteed whenever
\begin{equation*}
N^{-2\ell+}\delta^{\tfrac{1}{2}-2\ell-}N^{-3s} \frac{T}{\delta}\ll N^{-2s}.
\end{equation*}
So, it suffices that
\begin{equation*}
-2\ell+\frac{4s}{3}(\frac{1}{2}-2\ell)-3s-\frac{4s}{3}< -2s, \quad \textrm{i.e.}, \quad s>-6\ell/(5+8\ell).
\end{equation*}
Optimizing over the parameter $1/8<\ell<1/4$ (i.e., taking $\ell=1/4-$), we get $s>-3/14$. This completes the proof of theorem~\ref{global-theorem}.

\begin{remark}We take this opportunity to say that the $L^2\times L^2$ global well-posedness
of Schr\"odinger-Debye equation is \emph{not} proved in full details
in both papers~\cite{Corcho} and~\cite{AM}. Indeed, these papers
claim that the global well-posedness in $L^2\times L^2$ is an
immediate consequence of the conservation of the $L^2$-mass of $u$,
but they do not prevent a possible blow-up of $v$. However, it is
not hard to see that this can not occur in their context. In fact,
the $L^2$-norm of $v(t)$ can be controlled as follows:
\begin{equation*}
\begin{split}
\|v(t)\|_{L_x^2}&\leq e^{-t/\sigma}\|v_0\|_{L^2_x} +
\frac{1}{\sigma}\int_0^{t}e^{-(t-t')/\sigma}\|
|u(t')|^2\|_{L^2_x} dt'\\
&\leq e^{-t/\sigma}\|v_0\|_{L^2_x} +
\frac{1}{\sigma}\left\|e^{-(t-t')/\sigma}\right\|_{H_{t'}^{\tfrac{1}{2}-}}\| |u|^2 \|_{H_t^{-1/2+}L_x^2}\\
&\leq e^{-t/\sigma}\|v_0\|_{L^2_x} +
\frac{1}{\sigma}\left(\int_0^{t}e^{-2(t-t')/\sigma}dt'\right)^{1/2}
t^{1/4-}\| |u|^2\|_{H_t^{-1/4-}L_x^2}\\
&\leq e^{-t/\sigma}\|v_0\|_{L^2_x} +
C\left(1-e^{-2t/\sigma}\right)^{1/2}t^{3/4-}\|u\|_{X^{0,1/2+}}^{2}\\
&\leq e^{-t/\sigma}\|v_0\|_{L^2_x} +
Ct^{-1/2}\left(1-e^{-t/\sigma}\right)\cdot t^{3/4-}\|u_0\|_{L^2_x}^2\\
&\leq \max\{\|v_0\|_{L_x^2},C\|u_0\|_{L_x^2}^2\}
\end{split}
\end{equation*}
for any $0\leq t\leq 1$. Thus we have two scenarios:
\begin{itemize}
\item $\|v_0\|_{L_x^2}\leq C\|u_0\|_{L^2}^2$: in this situation, the
previous estimate implies that $\|v(t)\|_{L_x^2}\leq
C\|u_0\|_{L_x^2}^2$ for all $t$; since
$\|u(t)\|_{L^2}=\|u_0\|_{L^2}$ is a conserved quantity, there is no
blowup in this context;
\item $\|v_0\|_{L_x^2}\geq C\|u_0\|_{L^2}^{2}$: in this case, the previous estimate implies
$\|v(t)\|_{L^2}\leq \|v_0\|_{L^2}$ for all $t$ so that there is no
blowup occurring.
\end{itemize}
This completes the $L^2\times L^2$ global well-posedness arguments
of~\cite{Corcho} and~\cite{AM}.
\end{remark}


\section*{Acknowledgements}

The authors are thankful to professors F. Linares and B. Bid\'egaray for several
useful discussions concerning the Schr\"odinger-Debye equation, and the anonymous referee for pointing out a mistake in a previous version of this work.


\smallskip


\begin{thebibliography}{99}

\addcontentsline{toc}{chapter}{References}

\bibitem{ACH}{J. Angulo, A. Corcho and S. Hakkaev,}
{\it Well-posedness and stability of the periodic nonlinear waves interactions for the Benney system,}
{preprint 2008.}

\bibitem{AM}{A. Arbieto and C. Matheus,}
{\it On the periodic Schr\"odinger-Debye equation,}
{Comm. Pure and Applied Anal., {\bf 7} (2008), 699-713.}

\bibitem{Bekiranov1}{D. Bekiranov, T. Ogawa and G. Ponce,}
{\it Interaction equation for short and long dispersive waves,}
{J. Funct. Anal., {\bf 158} (1998), 357-388.}

\bibitem{Bidegaray}{B. Bid\'egaray}
{\it On the Cauchy problem for system ocurring in nonlinear optics,}
{Adv. Differential Equations, {\bf 3} (1998), 473-496.}

\bibitem{CKSTT2}{J. Colliander, M. Keel, G. Staffilani, H. Takaoka and T. Tao,}
{\it Multilinear estimates for periodic KdV equations, and applications,}
{J. Funct. Analysis, {\bf 211} (2004), 173-218.}

\bibitem{CKSTT3}{J. Colliander, M. Keel, G. Staffilani, H. Takaoka and T. Tao,}
{\it Global well-posedness for Schr\"odinger equations with derivative,}
{SIAM J. Math. Analysis, {\bf 33} (2001), 649-666.}

\bibitem{Corcho}{A. J. Corcho and F. Linares,}
{\it Well-Posedness for the Schr\"odinger-Debye Equation,}
{Contemporary Mathematics, {\bf 362} (2004), 113-131.}

\bibitem{Ginibre}{J. Ginibre, Y. Tsutsumi and G. Velo,}
{\it On the Cauchy Problem for the Zakharov system,}
{J. Funct. Anal., {\bf 151} (1997), 384-436.}

\bibitem{Grunrock}{A. Gr\"unrock,}
{\it An improved local well-posedness result for the modified KdV equation,}
{Int. Math. Res. Not., {\bf 61} (2004), 3287-3308.}

\bibitem{KPV}{C. Kenig, G. Ponce and L. Vega}
{\it Quadratic Forms for the 1-D Semilinear Schr\"odinger
Equation,} {Transactions of the American Mathematical Society,
{\bf 348} (1996), 3323-3353.}

\bibitem{T}{T. Tao}
{\it Nonlinear dispersive equations: local and global analysis,} {CBMS,
{\bf 106} (2006), 373 pp.}

\end{thebibliography}
\end{document}